\documentclass[10pt]{article}
\usepackage{amsmath}
\usepackage{amssymb}
\usepackage{amsthm}
\usepackage{mathrsfs}
\usepackage[compress,noadjust]{cite}
\numberwithin{equation}{section}

\begin{document}

\title{Endpoint estimates for commutators of intrinsic square functions in the Morrey type spaces}
\author{Hua Wang \footnote{E-mail address: wanghua@pku.edu.cn.}\\
\footnotesize{College of Mathematics and Econometrics, Hunan University, Changsha 410082, P. R. China}}
\date{}
\maketitle

\begin{abstract}
In this paper, the boundedness properties of commutators generated by $b$ and intrinsic square functions in the endpoint case are discussed, where $b\in BMO(\mathbb R^n)$. We first establish the weighted weak $L\log L$-type estimates for these commutator operators. Furthermore, we will prove endpoint estimates of commutators generated by $BMO(\mathbb R^n)$functions and intrinsic square functions in the weighted Morrey spaces $L^{1,\kappa}(w)$ for $0<\kappa<1$ and $w\in A_1$, and in the generalized Morrey spaces $L^{1,\Theta}$, where $\Theta$ is a growth function on $(0,+\infty)$ satisfying the doubling condition.
\end{abstract}

\section{Introduction and main results}

The intrinsic square functions were first introduced by Wilson in \cite{wilson1,wilson2}; they are defined as follows. For $0<\alpha\le1$, let ${\mathcal C}_\alpha$ be the family of functions $\varphi$ defined on $\mathbb R^n$ such that $\varphi$ has support containing in $\{x\in\mathbb R^n: |x|\le1\}$, $\int_{\mathbb R^n}\varphi(x)\,dx=0$, and for all $x, x'\in \mathbb R^n$,
\begin{equation*}
\big|\varphi(x)-\varphi(x')\big|\leq \big|x-x'\big|^{\alpha}.
\end{equation*}
For $(y,t)\in {\mathbb R}^{n+1}_{+}=\mathbb R^n\times(0,+\infty)$ and $f\in L^1_{{loc}}(\mathbb R^n)$, we set
\begin{equation}
A_\alpha(f)(y,t)=\sup_{\varphi\in{\mathcal C}_\alpha}\big|f*\varphi_t(y)\big|=\sup_{\varphi\in{\mathcal C}_\alpha}\bigg|\int_{\mathbb R^n}\varphi_t(y-z)f(z)\,dz\bigg|,
\end{equation}
where $\varphi_t(x)=t^{-n}\varphi(x/t)$. Then we define the intrinsic square function of $f$ (of order $\alpha$) by the formula
\begin{equation}
\mathcal S_{\alpha}(f)(x)=\left(\iint_{\Gamma(x)}\Big(A_\alpha(f)(y,t)\Big)^2\frac{dydt}{t^{n+1}}\right)^{1/2},
\end{equation}
where $\Gamma(x)$ denotes the usual cone of aperture one:
\begin{equation*}
\Gamma(x)=\big\{(y,t)\in{\mathbb R}^{n+1}_+:|x-y|<t\big\}.
\end{equation*}
Similarly, we can define a cone of aperture $\beta$ for any $\beta>0$:
\begin{equation*}
\Gamma_\beta(x)=\big\{(y,t)\in{\mathbb R}^{n+1}_+:|x-y|<\beta\cdot t\big\},
\end{equation*}
and the corresponding square function
\begin{equation}
\mathcal S_{\alpha,\beta}(f)(x)=\left(\iint_{\Gamma_\beta(x)}\Big(A_\alpha(f)(y,t)\Big)^2\frac{dydt}{t^{n+1}}\right)^{1/2}.
\end{equation}
The intrinsic Littlewood--Paley $\mathcal G$-function and the intrinsic $\mathcal G^*_\lambda$-function will be given respectively by
\begin{equation}
\mathcal G_\alpha(f)(x)=\left(\int_0^\infty\Big(A_\alpha(f)(x,t)\Big)^2\frac{dt}{t}\right)^{1/2}
\end{equation}
and
\begin{equation}
\mathcal G^*_{\lambda,\alpha}(f)(x)=\left(\iint_{{\mathbb R}^{n+1}_+}\left(\frac t{t+|x-y|}\right)^{\lambda n}\Big(A_\alpha(f)(y,t)\Big)^2\frac{dydt}{t^{n+1}}\right)^{1/2}, \quad \lambda>1.
\end{equation}

Let $b$ be a locally integrable function on $\mathbb R^n$, in this paper, we will consider the commutators generated by $b$ and intrinsic square functions, which are defined respectively by the following expressions in \cite{wang1}.
\begin{equation*}
\big[b,\mathcal S_\alpha\big](f)(x)=\left(\iint_{\Gamma(x)}\sup_{\varphi\in{\mathcal C}_\alpha}\bigg|\int_{\mathbb R^n}\big[b(x)-b(z)\big]\varphi_t(y-z)f(z)\,dz\bigg|^2\frac{dydt}{t^{n+1}}\right)^{1/2},
\end{equation*}
\begin{equation*}
\big[b,\mathcal G_\alpha\big](f)(x)=\left(\int_0^\infty\sup_{\varphi\in{\mathcal C}_\alpha}\bigg|\int_{\mathbb R^n}\big[b(x)-b(y)\big]\varphi_t(x-y)f(y)\,dy\bigg|^2\frac{dt}{t}\right)^{1/2},
\end{equation*}
and
\begin{equation*}
\begin{split}
&\big[b,\mathcal G^*_{\lambda,\alpha}\big](f)(x)\\
=&\left(\iint_{{\mathbb R}^{n+1}_+}\left(\frac t{t+|x-y|}\right)^{\lambda n}\sup_{\varphi\in{\mathcal C}_\alpha}\bigg|\int_{\mathbb R^n}\big[b(x)-b(z)\big]\varphi_t(y-z)f(z)\,dz\bigg|^2\frac{dydt}{t^{n+1}}\right)^{1/2}, \lambda>1.
\end{split}
\end{equation*}

On the other hand, the classical Morrey spaces $\mathcal L^{p,\lambda}$ were originally introduced by Morrey in \cite{morrey} to study the local behavior of solutions to second order elliptic partial differential equations. Since then, these spaces play an important role in studying the regularity of solutions to partial differential equations. For the boundedness of the Hardy--Littlewood maximal operator, the fractional integral operator and the Calder\'on--Zygmund singular integral operator on these spaces, we refer the reader to \cite{adams,chiarenza,peetre}. In \cite{mizuhara}, Mizuhara introduced the generalized Morrey space $L^{p,\Theta}$ which was later extended and studied by many authors (see \cite{guliyev1,guliyev2,guliyev3,lu,nakai}). In \cite{komori}, Komori and Shirai defined the weighted Morrey space $L^{p,\kappa}(w)$ which could be
viewed as an extension of weighted Lebesgue space, and then discussed the boundedness of the above classical operators in Harmonic Analysis on these weighted spaces. Recently, in \cite{wang1,wang2,wang3,wang4}, we have established the strong type and weak type estimates for intrinsic square functions and their commutators on $L^{p,\Theta}$ and $L^{p,\kappa}(w)$ with $1\leq p<\infty$.

In order to simplify the notations, for any given $\sigma>0$, we set
\begin{equation*}
\Phi\left(\frac{|f(x)|}{\sigma}\right)=\frac{|f(x)|}{\sigma}\cdot\left(1+\log^+\frac{|f(x)|}{\sigma}\right)
\end{equation*}
when $\Phi(t)=t\cdot(1+\log^+t)$. The main results of this paper can be stated as follows. For the endpoint estimates for these commutator operators $\big[b,\mathcal S_\alpha\big],\big[b,\mathcal G_{\alpha}\big]$ and $\big[b,\mathcal G^*_{\lambda,\alpha}\big]$ in the weighted Lebesgue space $L^1_w$,
when $b\in BMO(\mathbb R^n)$ and $w\in A_1$, we will obtain

\newtheorem{theorem}{Theorem}[section]

\begin{theorem}\label{mainthm:1}
Let $0<\alpha\leq1$, $w\in A_1$ and $b\in BMO(\mathbb R^n)$. Then for any given $\sigma>0$, there exists a constant $C>0$ independent of $f$ and $\sigma$ such that
\begin{equation*}
w\big(\big\{x\in\mathbb R^n:\big|[b,\mathcal S_\alpha](f)(x)\big|>\sigma\big\}\big)\leq C\int_{\mathbb R^n}\Phi\left(\frac{|f(x)|}{\sigma}\right)\cdot w(x)\,dx,
\end{equation*}
where $\Phi(t)=t(1+\log^+t)$.
\end{theorem}

\begin{theorem}\label{mainthm:2}
Let $0<\alpha\leq1$, $w\in A_1$ and $b\in BMO(\mathbb R^n)$. Then for any given $\sigma>0$, there exists a constant $C>0$ independent of $f$ and $\sigma$ such that
\begin{equation*}
w\big(\big\{x\in\mathbb R^n:\big|[b,\mathcal G_\alpha](f)(x)\big|>\sigma\big\}\big)\leq C\int_{\mathbb R^n}\Phi\left(\frac{|f(x)|}{\sigma}\right)\cdot w(x)\,dx,
\end{equation*}
where $\Phi(t)=t(1+\log^+t)$.
\end{theorem}

\begin{theorem}\label{mainthm:3}
Let $0<\alpha\leq1$, $w\in A_1$ and $b\in BMO(\mathbb R^n)$. If $\lambda>{(3n+2\alpha)}/n$, then for any given $\sigma>0$, there exists a constant $C>0$ independent of $f$ and $\sigma$ such that
\begin{equation*}
w\big(\big\{x\in\mathbb R^n:\big|[b,\mathcal G^*_{\lambda,\alpha}](f)(x)\big|>\sigma\big\}\big)\leq C\int_{\mathbb R^n}\Phi\left(\frac{|f(x)|}{\sigma}\right)\cdot w(x)\,dx,
\end{equation*}
where $\Phi(t)=t(1+\log^+t)$.
\end{theorem}

In particular, if we take $w$ to be a constant function, then we immediately
get the following:

\newtheorem{corollary}[theorem]{Corollary}

\begin{corollary}\label{maincor:1}
Let $0<\alpha\leq1$ and $b\in BMO(\mathbb R^n)$. Then for any given $\sigma>0$, there exists a constant $C>0$ independent of $f$ and $\sigma$ such that
\begin{equation*}
\big|\big\{x\in\mathbb R^n:\big|[b,\mathcal S_\alpha](f)(x)\big|>\sigma\big\}\big|\leq C\int_{\mathbb R^n}\Phi\left(\frac{|f(x)|}{\sigma}\right)dx,
\end{equation*}
where $\Phi(t)=t(1+\log^+t)$.
\end{corollary}

\begin{corollary}\label{maincor:2}
Let $0<\alpha\leq1$ and $b\in BMO(\mathbb R^n)$. Then for any given $\sigma>0$, there exists a constant $C>0$ independent of $f$ and $\sigma$ such that
\begin{equation*}
\big|\big\{x\in\mathbb R^n:\big|[b,\mathcal G_\alpha](f)(x)\big|>\sigma\big\}\big|\leq C\int_{\mathbb R^n}\Phi\left(\frac{|f(x)|}{\sigma}\right)dx,
\end{equation*}
where $\Phi(t)=t(1+\log^+t)$.
\end{corollary}

\begin{corollary}\label{maincor:3}
Let $0<\alpha\leq1$ and $b\in BMO(\mathbb R^n)$. If $\lambda>{(3n+2\alpha)}/n$, then for any given $\sigma>0$, there exists a constant $C>0$ independent of $f$ and $\sigma$ such that
\begin{equation*}
\big|\big\{x\in\mathbb R^n:\big|[b,\mathcal G^*_{\lambda,\alpha}](f)(x)\big|>\sigma\big\}\big|\leq C\int_{\mathbb R^n}\Phi\left(\frac{|f(x)|}{\sigma}\right)dx,
\end{equation*}
where $\Phi(t)=t(1+\log^+t)$.
\end{corollary}

For the endpoint estimates of commutators generated by $BMO(\mathbb R^n)$ functions and intrinsic square functions in the weighted Morrey spaces $L^{1,\kappa}(w)$ for all $0<\kappa<1$ and $w\in A_1$, we will prove

\begin{theorem}\label{mainthm:4}
Let $0<\alpha\leq1$, $0<\kappa<1$, $w\in A_1$ and $b\in BMO(\mathbb R^n)$. Then for any given $\sigma>0$ and any ball $B$, there exists a constant $C>0$ independent of $f$, $B$ and $\sigma$ such that
\begin{equation*}
\frac{1}{w(B)^\kappa}\cdot w\big(\big\{x\in B:\big|[b,\mathcal S_\alpha](f)(x)\big|>\sigma\big\}\big)\leq C\cdot\sup_B\frac{1}{w(B)^\kappa}
\int_{B}\Phi\left(\frac{|f(x)|}{\sigma}\right)\cdot w(x)dx,
\end{equation*}
where $\Phi(t)=t(1+\log^+t)$.
\end{theorem}

\begin{theorem}\label{mainthm:5}
Let $0<\alpha\leq1$, $0<\kappa<1$, $w\in A_1$ and $b\in BMO(\mathbb R^n)$. Then for any given $\sigma>0$ and any ball $B$, there exists a constant $C>0$ independent of $f$, $B$ and $\sigma$ such that
\begin{equation*}
\frac{1}{w(B)^\kappa}\cdot w\big(\big\{x\in B:\big|[b,\mathcal G_\alpha](f)(x)\big|>\sigma\big\}\big)\leq C\cdot\sup_B\frac{1}{w(B)^\kappa}
\int_{B}\Phi\left(\frac{|f(x)|}{\sigma}\right)\cdot w(x)dx,
\end{equation*}
where $\Phi(t)=t(1+\log^+t)$.
\end{theorem}

\begin{theorem}\label{mainthm:6}
Let $0<\alpha\leq1$, $0<\kappa<1$, $w\in A_1$ and $b\in BMO(\mathbb R^n)$. If $\lambda>{(3n+2\alpha)}/n$, then for any given $\sigma>0$ and any ball $B$, there exists a constant $C>0$ independent of $f$, $B$ and $\sigma$ such that
\begin{equation*}
\frac{1}{w(B)^\kappa}\cdot w\big(\big\{x\in B:\big|[b,\mathcal G^*_{\lambda,\alpha}](f)(x)\big|>\sigma\big\}\big)\leq C\cdot\sup_B\frac{1}{w(B)^\kappa}
\int_{B}\Phi\left(\frac{|f(x)|}{\sigma}\right)\cdot w(x)dx,
\end{equation*}
where $\Phi(t)=t(1+\log^+t)$.
\end{theorem}

For the endpoint estimates of commutators generated by $BMO(\mathbb R^n)$ functions and intrinsic square functions in the generalized Morrey spaces $L^{1,\Theta}$ when $\Theta$ satisfies the doubling condition, we will show that

\begin{theorem}\label{mainthm:7}
Let $0<\alpha\leq1$ and $b\in BMO(\mathbb R^n)$. Suppose that $\Theta$ satisfies $(\ref{doubling})$ and $1\le D(\Theta)<2^n$, then for any given $\sigma>0$ and any ball $B(x_0,r)$, there exists a constant $C>0$ independent of $f$, $B(x_0,r)$ and $\sigma$ such that
\begin{equation*}
\frac{1}{\Theta(r)}\cdot\big|\big\{x\in B(x_0,r):\big|[b,\mathcal S_\alpha](f)(x)\big|>\sigma\big\}\big|\leq C\cdot\sup_{r>0}\frac{1}{\Theta(r)}
\int_{B(x_0,r)}\Phi\left(\frac{|f(x)|}{\sigma}\right)dx,
\end{equation*}
where $\Phi(t)=t(1+\log^+t)$.
\end{theorem}

\begin{theorem}\label{mainthm:8}
Let $0<\alpha\leq1$ and $b\in BMO(\mathbb R^n)$. Suppose that $\Theta$ satisfies $(\ref{doubling})$ and $1\le D(\Theta)<2^n$, then for any given $\sigma>0$ and any ball $B(x_0,r)$, there exists a constant $C>0$ independent of $f$, $B(x_0,r)$ and $\sigma$ such that
\begin{equation*}
\frac{1}{\Theta(r)}\cdot\big|\big\{x\in B(x_0,r):\big|[b,\mathcal G_\alpha](f)(x)\big|>\sigma\big\}\big|\leq C\cdot\sup_{r>0}\frac{1}{\Theta(r)}
\int_{B(x_0,r)}\Phi\left(\frac{|f(x)|}{\sigma}\right)dx,
\end{equation*}
where $\Phi(t)=t(1+\log^+t)$.
\end{theorem}

\begin{theorem}\label{mainthm:9}
Let $0<\alpha\leq1$ and $b\in BMO(\mathbb R^n)$. Suppose that $\Theta$ satisfies $(\ref{doubling})$, $1\le D(\Theta)<2^n$ and $\lambda>{(3n+2\alpha)}/n$, then for any given $\sigma>0$ and any ball $B(x_0,r)$, there exists a constant $C>0$ independent of $f$, $B(x_0,r)$ and $\sigma$ such that
\begin{equation*}
\frac{1}{\Theta(r)}\cdot\big|\big\{x\in B(x_0,r):\big|[b,\mathcal G^*_{\lambda,\alpha}](f)(x)\big|>\sigma\big\}\big|\leq C\cdot\sup_{r>0}\frac{1}{\Theta(r)}
\int_{B(x_0,r)}\Phi\left(\frac{|f(x)|}{\sigma}\right)dx,
\end{equation*}
where $\Phi(t)=t(1+\log^+t)$.
\end{theorem}

\section{Notations and preliminaries}

A weight $w$ will always mean a positive function which is locally integrable on $\mathbb R^n$, $B=B(x_0,r_B)$ denotes the open ball with the center $x_0$ and radius $r_B$. For $1<p<\infty$, a weight function $w$ is said to belong to the Muckenhoupt's class $A_p$, if there is a constant $C>0$ such that for every ball $B\subseteq \mathbb R^n$(see \cite{garcia,muckenhoupt}),
\begin{equation*}
\left(\frac1{|B|}\int_B w(x)\,dx\right)\left(\frac1{|B|}\int_B w(x)^{-1/{(p-1)}}\,dx\right)^{p-1}\le C.
\end{equation*}
For the case $p=1$, $w\in A_1$, if there is a constant $C>0$ such that for every ball $B\subseteq \mathbb R^n$,
\begin{equation*}
\frac1{|B|}\int_B w(x)\,dx\le C\cdot\underset{x\in B}{\mbox{ess\,inf}}\;w(x).
\end{equation*}
We also define $A_\infty=\cup_{1\leq p<\infty}A_p$. It is well known that if $w\in A_p$ with $1\leq p<\infty$, then for any ball $B$, there exists an absolute constant $C>0$ such that
\begin{equation}\label{weights}
w(2B)\le C\,w(B).
\end{equation}
In general, for $w\in A_1$ and any $j\in\mathbb Z_+$, there exists an absolute constant $C>0$ such that (see \cite{garcia})
\begin{equation}\label{general weights}
w\big(2^j B\big)\le C\cdot 2^{jn}w(B).
\end{equation}
Moreover, if $w\in A_\infty$, then for all balls $B$ and all measurable subsets $E$ of $B$, there exists a number $\delta>0$ independent of $E$ and $B$ such that (see \cite{garcia})
\begin{equation}\label{compare}
\frac{w(E)}{w(B)}\le C\left(\frac{|E|}{|B|}\right)^\delta.
\end{equation}
A weight function $w$ is said to belong to the reverse H\"{o}lder class $RH_r$, if there exist two constants $r>1$ and
$C>0$ such that the following reverse H\"{o}lder inequality holds for every
ball $B\subseteq \mathbb R^n$.
\begin{equation*}
\left(\frac{1}{|B|}\int_B w(x)^r\,dx\right)^{1/r}\le C\left(\frac{1}{|B|}\int_B w(x)\,dx\right).
\end{equation*}
Given a ball $B$ and $\lambda>0$, $\lambda B$ denotes the ball with the same center as $B$ whose radius is $\lambda$ times that of $B$. For a given weight function $w$ and a measurable set $E$, we also denote the Lebesgue measure of $E$ by $|E|$ and the weighted measure of $E$ by $w(E)$, where $w(E)=\int_E w(x)\,dx$. Equivalently, we could define the above notions with cubes
instead of balls. Hence we shall use these two different definitions appropriate to calculations.

Given a weight function $w$ on $\mathbb R^n$, for $1\leq p<\infty$, the weighted Lebesgue space $L^p_w(\mathbb R^n)$ is defined as the set of all functions $f$ such that
\begin{equation*}
\big\|f\big\|_{L^p_w}=\bigg(\int_{\mathbb R^n}|f(x)|^pw(x)\,dx\bigg)^{1/p}<\infty.
\end{equation*}
In particular, when $w$ equals to a constant function, we will denote $L^p_w(\mathbb R^n)$ simply by $L^p(\mathbb R^n)$.

Let $0<\kappa<1$ and $w$ be a weight function on $\mathbb R^n$. Then the weighted Morrey space $L^{1,\kappa}(w)$ is defined by (see \cite{komori})
\begin{equation*}
L^{1,\kappa}(w)=\left\{f\in L^1_{loc}(w):\big\|f\big\|_{L^{1,\kappa}(w)}=\sup_B\frac{1}{w(B)^{\kappa}}\int_B|f(x)|w(x)\,dx<\infty\right\},
\end{equation*}
where the supremum is taken over all balls $B$ in $\mathbb R^n$.

Let $\Theta=\Theta(r)$, $r>0$, be a growth function, that is, a positive increasing function in $(0,+\infty)$ and satisfy the following doubling condition:
\begin{equation}\label{doubling}
\Theta(2r)\leq D\cdot\Theta(r), \quad \mbox{for all }\,r>0,
\end{equation}
where $D=D(\Theta)\ge1$ is a doubling constant independent of $r$.
The generalized Morrey space $L^{1,\Theta}(\mathbb R^n)$ is defined as the set of all locally integrable functions $f$ for which (see \cite{mizuhara})
\begin{equation*}
\sup_{r>0;B(x_0,r)}\frac{1}{\Phi(r)}\int_{B(x_0,r)}|f(x)|\,dx<\infty,
\end{equation*}
where $B(x_0,r)=\{x\in\mathbb R^n:|x-x_0|<r\}$ is the open ball centered at $x_0$ and with radius $r>0$.

We next recall some basic definitions and facts about Orlicz spaces needed for the proof of the main results. For more information on the subject, one can see \cite{rao}. A function $\Phi$ is called a Young function if it is continuous, nonnegative, convex and strictly increasing on $[0,+\infty)$ with $\Phi(0)=0$ and $\Phi(t)\to +\infty$ as $t\to +\infty$. We define the $\Phi$-average of a function $f$ over a ball $B$ by means of the following Luxemburg norm:
\begin{equation*}
\big\|f\big\|_{\Phi,B}=\inf\left\{\sigma>0:\frac{1}{|B|}\int_B\Phi\left(\frac{|f(x)|}{\sigma}\right)dx\leq1\right\}.
\end{equation*}
An equivalent norm that is often useful in calculations is as follows(see \cite{rao,perez1}):
\begin{equation}\label{equiv norm}
\big\|f\big\|_{\Phi,B}\leq \inf_{\eta>0}\left\{\eta+\frac{\eta}{|B|}\int_B\Phi\left(\frac{|f(x)|}{\eta}\right)dx\right\}\leq 2\big\|f\big\|_{\Phi,B}.
\end{equation}
Given a Young function $\Phi$, we use $\bar\Phi$ to denote the complementary Young function associated to $\Phi$. Then the following generalized H\"older's inequality holds for any given ball $B$ (see \cite{perez1,perez2}).
\begin{equation*}
\frac{1}{|B|}\int_B|f(x)\cdot g(x)|dx\leq 2\big\|f\big\|_{\Phi,B}\big\|g\big\|_{\bar\Phi,B}.
\end{equation*}
In order to deal with the weighted case, for $w\in A_\infty$, we also need to define the weighted $\Phi$-average of a function $f$ over a ball $B$ by means of the weighted Luxemburg norm:
\begin{equation*}
\big\|f\big\|_{\Phi(w),B}=\inf\left\{\sigma>0:\frac{1}{w(B)}\int_B\Phi\left(\frac{|f(x)|}{\sigma}\right)w(x)\,dx\leq1\right\}.
\end{equation*}
It can be shown that for $w\in A_\infty$(see \cite{rao,zhang}),
\begin{equation}\label{equiv norm with weight}
\big\|f\big\|_{\Phi(w),B}\approx \inf_{\eta>0}\left\{\eta+\frac{\eta}{w(B)}\int_B\Phi\left(\frac{|f(x)|}{\eta}\right)w(x)\,dx\right\},
\end{equation}
and
\begin{equation*}
\frac{1}{w(B)}\int_B|f(x)g(x)|w(x)\,dx\leq C\big\|f\big\|_{\Phi(w),B}\big\|g\big\|_{\bar\Phi(w),B}.
\end{equation*}
The young function that we are going to use is $\Phi(t)=t(1+\log^+t)$ with its complementary Young function $\bar\Phi(t)\approx\exp(t)$. In the present situation, we denote
\begin{equation*}
\big\|f\big\|_{L\log L,B}=\big\|f\big\|_{\Phi,B}, \qquad \big\|g\big\|_{\exp L,B}=\big\|g\big\|_{\bar\Phi,B};
\end{equation*}
and
\begin{equation*}
\big\|f\big\|_{L\log L(w),B}=\big\|f\big\|_{\Phi(w),B}, \qquad \big\|g\big\|_{\exp L(w),B}=\big\|g\big\|_{\bar\Phi(w),B}.
\end{equation*}
By the (weighted) generalized H\"older's inequality, we have (see \cite{perez1,zhang})
\begin{equation}\label{holder}
\frac{1}{|B|}\int_B|f(x)\cdot g(x)|dx\leq 2\big\|f\big\|_{L\log L,B}\big\|g\big\|_{\exp L,B},
\end{equation}
and
\begin{equation}\label{weighted holder}
\frac{1}{w(B)}\int_B|f(x)g(x)|w(x)\,dx\leq C\big\|f\big\|_{L\log L(w),B}\big\|g\big\|_{\exp L(w),B}.
\end{equation}

Let us now recall the definition of the space of $BMO(\mathbb R^n)$ (Bounded Mean Oscillation) (see \cite{duoand,john}).
A locally integrable function $b$ is said to be in $BMO(\mathbb R^n)$, if
\begin{equation*}
\|b\|_*=\sup_{B}\frac{1}{|B|}\int_B|b(x)-b_B|\,dx<\infty,
\end{equation*}
where $b_B$ stands for the average of $b$ on $B$, i.e., $b_B=\frac{1}{|B|}\int_B b(y)\,dy$ and the supremum is taken
over all balls $B$ in $\mathbb R^n$. Modulo constants, the space $BMO(\mathbb R^n)$ is a Banach space with respect to the norm $\|\cdot\|_*$.
By the John--Nirenberg's inequality, it is not difficult to see that for any given ball $B$ (see \cite{perez1,perez2})
\begin{equation}\label{exp}
\big\|b-b_B\big\|_{\exp L,B}\leq C\|b\|_*.
\end{equation}
Furthermore, we can also prove that for any $w\in A_\infty$ and any given ball $B$ (see \cite{zhang}),
\begin{equation}\label{weighted exp}
\big\|b-b_B\big\|_{\exp L(w),B}\leq C\|b\|_*.
\end{equation}

Throughout this paper, the letter $C$ always denotes a positive constant independent of the main parameters involved, but it may be different from line to line. By $A\approx B$, we mean that there exists a constant $C>1$ such that $\frac1C\le\frac AB\le C$.

\section{Proofs of Theorems \ref{mainthm:1} and \ref{mainthm:2}}

Given a real-valued function $b\in BMO(\mathbb R^n)$, we shall follow the idea developed in \cite{alvarez,ding} and denote $F(\xi)=e^{\xi[b(x)-b(z)]}$, $\xi\in\mathbb C$. Then by the analyticity of $F(\xi)$ on $\mathbb C$ and the Cauchy integral formula, we get
\begin{equation*}
\begin{split}
b(x)-b(z)=F'(0)&=\frac{1}{2\pi i}\int_{|\xi|=1}\frac{F(\xi)}{\xi^2}\,d\xi\\
&=\frac{1}{2\pi}\int_0^{2\pi}e^{e^{i\theta}[b(x)-b(z)]}\cdot e^{-i\theta}\,d\theta.
\end{split}
\end{equation*}
Thus, for any $\varphi\in{\mathcal C}_\alpha$, $0<\alpha\le1$, we obtain
\begin{equation*}
\begin{split}
&\bigg|\int_{\mathbb R^n}\big[b(x)-b(z)\big]\varphi_t(y-z)f(z)\,dz\bigg|\\
=&
\bigg|\frac{1}{2\pi}\int_0^{2\pi}\bigg(\int_{\mathbb R^n}\varphi_t(y-z)e^{-e^{i\theta}b(z)}f(z)\,dz\bigg)
e^{e^{i\theta}b(x)}\cdot e^{-i\theta}\,d\theta\bigg|\\
\leq&\frac{1}{2\pi}\int_0^{2\pi}\sup_{\varphi\in{\mathcal C}_\alpha}\bigg|\int_{\mathbb R^n}\varphi_t(y-z)e^{-e^{i\theta}b(z)}f(z)\,dz\bigg|e^{\cos\theta\cdot b(x)}\,d\theta\\
\leq&\frac{1}{2\pi}\int_0^{2\pi}A_\alpha\big(e^{-e^{i\theta}b}\cdot f\big)(y,t)\cdot e^{\cos\theta\cdot b(x)}\,d\theta.
\end{split}
\end{equation*}
So we have
\begin{equation*}
\big|\big[b,\mathcal S_\alpha\big](f)(x)\big|\le\frac{1}{2\pi}\int_0^{2\pi}
\mathcal S_\alpha\big(e^{-e^{i\theta}b}\cdot f\big)(x)\cdot e^{\cos\theta\cdot b(x)}\,d\theta,
\end{equation*}
and
\begin{equation*}
\big|\big[b,\mathcal G_\alpha\big](f)(x)\big|\le\frac{1}{2\pi}\int_0^{2\pi}
\mathcal G_\alpha\big(e^{-e^{i\theta}b}\cdot f\big)(x)\cdot e^{\cos\theta\cdot b(x)}\,d\theta.
\end{equation*}

Then, by the $L^p_w$-boundedness of intrinsic square functions (see \cite{wilson2}), and using the same arguments as in \cite{ding}, we can also show the following:

\begin{theorem}\label{commutator thm}
Let $0<\alpha\le1$, $1<p<\infty$ and $w\in A_p$. Then the commutators $\big[b,\mathcal S_\alpha\big]$ and $\big[b,\mathcal G_{\alpha}\big]$ are all bounded from $L^p_w(\mathbb R^n)$ into itself whenever $b\in BMO(\mathbb R^n)$.
\end{theorem}

We are now ready to give the proofs of Theorems $\ref{mainthm:1}$ and $\ref{mainthm:2}$, which are based on the Calder\'on--Zygmund decomposition.
\begin{proof}[Proofs of Theorems $\ref{mainthm:1}$ and $\ref{mainthm:2}$]
We will only give the proof of Theorem $\ref{mainthm:1}$ here, since the proof of Theorem $\ref{mainthm:2}$ is similar and easier. Inspired by the work in \cite{perez2,perez3,zhang}, for any fixed $\sigma>0$, we apply the Calder\'on--Zygmund decomposition of $f$ at height $\sigma$ to obtain a sequence of disjoint non-overlapping dyadic cubes $\{Q_i\}$ such that the following property holds (see \cite{stein})
\begin{equation}\label{decomposition}
\sigma<\frac{1}{|Q_i|}\int_{Q_i}|f(y)|\,dy< 2^n\cdot\sigma,
\end{equation}
where $Q_i=Q(c_i,\ell_i)$ denotes the cube centered at $c_i$ with side length $\ell_i$ and all cubes are
assumed to have their sides parallel to the coordinate axes. Setting $E=\bigcup_i Q_i$. Now we define two functions $g$ and $h$ as follows:
\begin{equation*}
g(x)=
\begin{cases}
f(x) &  \mbox{if}\;\; x\in E^c,\\
\frac{1}{|Q_i|}\int_{Q_i}|f(y)|\,dy    &  \mbox{if}\;\; x\in Q_i,
\end{cases}
\end{equation*}
and
\begin{equation*}
h(x)=f(x)-g(x)=\sum_i h_i(x),
\end{equation*}
where $h_i(x)=h(x)\chi_{Q_i}(x)$. Then we have
\begin{equation}\label{pointwise estimate g}
|g(x)|\le C\cdot\sigma, \quad \mbox{a.e. }\, x\in\mathbb R^n,
\end{equation}
and
\begin{equation}\label{f=g+h}
f(x)=g(x)+h(x).
\end{equation}
Obviously, supp\,$h_i\subseteq Q_i$, $\int_{Q_i}h_i(x)\,dx=0$ and $\|h_i\|_{L^1}\le 2\int_{Q_i}|f(x)|\,dx$ by the above decomposition. Since $\big|\big[b,\mathcal S_\alpha\big](f)(x)\big|\leq\big|\big[b,\mathcal S_\alpha\big](g)(x)\big|+\big|\big[b,\mathcal S_\alpha\big](h)(x)\big|$ by (\ref{f=g+h}), then we can write
\begin{equation*}
\begin{split}
&w\big(\big\{x\in\mathbb R^n:\big|\big[b,\mathcal S_\alpha\big](f)(x)\big|>\sigma\big\}\big)\\
\leq&
w\big(\big\{x\in\mathbb R^n:\big|\big[b,\mathcal S_\alpha\big](g)(x)\big|>\sigma/2\big\}\big)
+w\big(\big\{x\in\mathbb R^n:\big|\big[b,\mathcal S_\alpha\big](h)(x)\big|>\sigma/2\big\}\big)\\
:=&I_1+I_2.
\end{split}
\end{equation*}
Observe that $w\in A_1\subset A_2$. Applying Chebyshev's inequality and Theorem \ref{commutator thm}, we obtain
\begin{equation*}
I_1\leq \frac{4}{\sigma^2}\cdot\Big\|\big[b,\mathcal S_\alpha\big](g)\Big\|^2_{L^2_w}\leq \frac{C}{\sigma^2}\cdot\big\|g\big\|^2_{L^2_w}.
\end{equation*}
Moreover, by the inequality (\ref{pointwise estimate g}) and the $A_1$ condition, we deduce that
\begin{align}\label{g}
\big\|g\big\|^2_{L^2_w}&\le C\cdot\sigma\int_{\mathbb R^n}|g(x)|w(x)\,dx\notag\\
&\le C\cdot\sigma\left(\int_{E^c}|f(x)|w(x)\,dx+\int_{\bigcup_i Q_i}|g(x)|w(x)\,dx\right)\notag\\
&\le C\cdot\sigma\left(\int_{\mathbb R^n}|f(x)|w(x)\,dx+\sum_i\frac{w(Q_i)}{|Q_i|}\int_{Q_i}|f(y)|\,dy\right)\notag\\
&\le C\cdot\sigma\left(\int_{\mathbb R^n}|f(x)|w(x)\,dx+\sum_i\underset{y\in Q_i}{\mbox{ess\,inf}}\,w(y)\int_{Q_i}|f(y)|\,dy\right)\notag\\
&\le C\cdot\sigma\left(\int_{\mathbb R^n}|f(x)|w(x)\,dx+\int_{\bigcup_i Q_i}|f(y)|w(y)\,dy\right)\notag\\
&\le C\cdot\sigma\int_{\mathbb R^n}|f(x)|w(x)\,dx.
\end{align}
So we have
\begin{equation*}
I_1\leq C\int_{\mathbb R^n}\frac{|f(x)|}{\sigma}\cdot w(x)\,dx\leq C\int_{\mathbb R^n}\Phi\left(\frac{|f(x)|}{\sigma}\right)\cdot w(x)\,dx.
\end{equation*}
To deal with the other term $I_2$, let $Q_i^*=2\sqrt n Q_i$ be the cube concentric with $Q_i$ such that $\ell(Q_i^*)=(2\sqrt n)\ell(Q_i)$. Then we can further decompose $I_2$ as follows.
\begin{equation*}
\begin{split}
I_2\le&\,w\Big(\Big\{x\in \bigcup_i Q_i^*:\Big|\big[b,\mathcal S_\alpha\big](h)(x)\Big|>\sigma/2\Big\}\Big)\\
&+w\Big(\Big\{x\notin \bigcup_i Q_i^*:\Big|\big[b,\mathcal S_\alpha\big](h)(x)\Big|>\sigma/2\Big\}\Big)\\
:=&\,I_3+I_4.
\end{split}
\end{equation*}
Since $w\in A_1$, then by the inequality (\ref{weights}), we can get
\begin{equation*}
I_3\leq\sum_i w\big(Q_i^*\big)\le C\sum_i w(Q_i).
\end{equation*}
Furthermore, it follows from the inequality (\ref{decomposition}) and the $A_1$ condition that
\begin{equation*}
\begin{split}
I_3&\leq C\sum_i\frac{\,1\,}{\sigma}\cdot\underset{y\in Q_i}{\mbox{ess\,inf}}\,w(y)\int_{Q_i}|f(y)|\,dy\\
&\leq \frac{C}{\sigma}\sum_i\int_{Q_i}|f(y)|w(y)\,dy\leq \frac{C}{\sigma}\int_{\bigcup_i Q_i}|f(y)|w(y)\,dy\\
&\leq C\int_{\mathbb R^n}\frac{|f(y)|}{\sigma}\cdot w(y)\,dy\leq C\int_{\mathbb R^n}\Phi\left(\frac{|f(y)|}{\sigma}\right)\cdot w(y)\,dy.
\end{split}
\end{equation*}
For any given $x\in \mathbb R^n$ and $(y,t)\in\Gamma(x)$, we have
\begin{align}\label{commutator estimate}
\sup_{\varphi\in{\mathcal C}_\alpha}\bigg|\int_{\mathbb R^n}\big[b(x)-b(z)\big]\varphi_t(y-z)h_i(z)\,dz\bigg|&\leq
\big|b(x)-b_{Q_i}\big|\cdot\sup_{\varphi\in{\mathcal C}_\alpha}\bigg|\int_{\mathbb R^n}\varphi_t(y-z)h_i(z)\,dz\bigg|\notag\\
&+\sup_{\varphi\in{\mathcal C}_\alpha}\bigg|\int_{\mathbb R^n}\big[b(z)-b_{Q_i}\big]\varphi_t(y-z)h_i(z)\,dz\bigg|.
\end{align}
Hence
\begin{equation*}
\begin{split}
\big|\big[b,\mathcal S_\alpha\big](h)(x)\big|&\leq\sum_i\big|b(x)-b_{Q_i}\big|\cdot\mathcal S_\alpha(h_i)(x)\\
&+\left(\iint_{\Gamma(x)}\sup_{\varphi\in{\mathcal C}_\alpha}\bigg|\int_{\mathbb R^n}\big[b(z)-b_{Q_i}\big]\varphi_t(y-z)\cdot\sum_i h_i(z)\,dz\bigg|^2\frac{dydt}{t^{n+1}}\right)^{1/2}\\
&=\sum_i\big|b(x)-b_{Q_i}\big|\cdot\mathcal S_\alpha(h_i)(x)+\mathcal S_\alpha\bigg(\sum_i[b-b_{Q_i}]h_i\bigg)(x).
\end{split}
\end{equation*}
Then we can write
\begin{equation*}
\begin{split}
I_4\leq & w\bigg(\bigg\{x\notin \bigcup_i Q_i^*:\sum_i\big|b(x)-b_{Q_i}\big|\cdot\mathcal S_\alpha(h_i)(x)>\sigma/4\bigg\}\bigg)\\
&+ w\bigg(\bigg\{x\notin \bigcup_i Q_i^*:\mathcal S_\alpha\bigg(\sum_i[b-b_{Q_i}]h_i\bigg)(x)>\sigma/4\bigg\}\bigg)\\
:=&I_5+I_6.
\end{split}
\end{equation*}
It follows directly from the Chebyshev's inequality that
\begin{equation*}
\begin{split}
I_5&\leq\frac{\,4\,}{\sigma}\int_{\mathbb R^n\backslash\bigcup_i Q_i^*}\bigg|\sum_i\big|b(x)-b_{Q_i}\big|\cdot\mathcal S_\alpha(h_i)(x)\bigg|w(x)\,dx\\
&\leq\frac{\,4\,}{\sigma}\sum_i\left(\int_{(Q_i^*)^c}\big|b(x)-b_{Q_i}\big|\cdot\mathcal S_\alpha(h_i)(x)\,w(x)\,dx\right).\\
\end{split}
\end{equation*}
Denote the center of $Q_i$ by $c_i$. For any $\varphi\in{\mathcal C}_\alpha$, $0<\alpha\le1$, by the cancellation condition of $h_i$, we obtain that for any $(y,t)\in\Gamma(x)$,
\begin{align}\label{kernel estimate1}
\big|(h_i*\varphi_t)(y)\big|&=\left|\int_{Q_i}\big[\varphi_t(y-z)-\varphi_t(y-c_i)\big]h_i(z)\,dz\right|\notag\\
&\le\int_{Q_i\cap\{z:|z-y|\le t\}}\frac{|z-c_i|^\alpha}{t^{n+\alpha}}|h_i(z)|\,dz\notag\\
&\le C\cdot\frac{\ell(Q_i)^{\alpha}}{t^{n+\alpha}}\int_{Q_i\cap\{z:|z-y|\le t\}}|h_i(z)|\,dz.
\end{align}
In addition, for any $z\in Q_i$ and $x\in (Q^*_i)^c$, we have $|z-c_i|<\frac{|x-c_i|}{2}$. Thus, for all $(y,t)\in\Gamma(x)$ and $|z-y|\le t$ with $z\in Q_i$, we can see that
\begin{equation}\label{2t}
t+t\ge|x-y|+|y-z|\ge|x-z|\ge|x-c_i|-|z-c_i|\ge\frac{|x-c_i|}{2}.
\end{equation}
Hence, for any $x\in (Q^*_i)^c$, by using the above inequalities (\ref{kernel estimate1}) and (\ref{2t}), we obtain
\begin{equation*}
\begin{split}
\big|\mathcal S_{\alpha}(h_i)(x)\big|&=\left(\iint_{\Gamma(x)}\bigg(\sup_{\varphi\in{\mathcal C}_\alpha}\big|(\varphi_t*{h_i})(y)\big|\bigg)^2\frac{dydt}{t^{n+1}}\right)^{1/2}\\
&\leq C\cdot\ell(Q_i)^{\alpha}\bigg(\int_{Q_i}|h_i(z)|\,dz\bigg)\left(\int_{\frac{|x-c_i|}{4}}^\infty
\int_{|y-x|<t}\frac{dydt}{t^{2(n+\alpha)+n+1}}\right)^{1/2}\\
&\leq C\cdot\ell(Q_i)^{\alpha}\bigg(\int_{Q_i}|h_i(z)|\,dz\bigg)
\left(\int_{\frac{|x-c_i|}{4}}^\infty\frac{dt}{t^{2(n+\alpha)+1}}\right)^{1/2}\\
&\leq C\cdot\frac{\ell(Q_i)^{\alpha}}{|x-c_i|^{n+\alpha}}\bigg(\int_{Q_i}|f(z)|\,dz\bigg).
\end{split}
\end{equation*}
Since $Q_i^*=2\sqrt n Q_i\supset 2Q_i$, then $(Q_i^*)^c\subset (2Q_i)^c$. This fact together with the above pointwise estimate yields
\begin{equation*}
\begin{split}
I_5&\leq\frac{C}{\sigma}\sum_i
\left(\ell(Q_i)^{\alpha}\int_{Q_i}|f(z)|\,dz\times\int_{(Q_i^*)^c}\big|b(x)-b_{Q_i}\big|\cdot\frac{w(x)}{|x-c_i|^{n+\alpha}}dx\right)\\
&\leq\frac{C}{\sigma}\sum_i
\left(\ell(Q_i)^{\alpha}\int_{Q_i}|f(z)|\,dz\times\int_{(2Q_i)^c}\big|b(x)-b_{Q_i}\big|\cdot\frac{w(x)}{|x-c_i|^{n+\alpha}}dx\right)\\
&\leq\frac{C}{\sigma}\sum_i
\left(\ell(Q_i)^{\alpha}\int_{Q_i}|f(z)|\,dz\times\sum_{j=1}^\infty\int_{2^{j+1}Q_i\backslash 2^{j}Q_i}\big|b(x)-b_{2^{j+1}Q_i}\big|\cdot\frac{w(x)}{|x-c_i|^{n+\alpha}}dx\right)\\
&+\frac{C}{\sigma}\sum_i
\left(\ell(Q_i)^{\alpha}\int_{Q_i}|f(z)|\,dz\times\sum_{j=1}^\infty\int_{2^{j+1}Q_i\backslash 2^{j}Q_i}\big|b_{2^{j+1}Q_i}-b_{Q_i}\big|\cdot\frac{w(x)}{|x-c_i|^{n+\alpha}}dx\right)\\
&:=\mbox{\upshape I+II}.
\end{split}
\end{equation*}
For the term I,
\begin{equation*}
\begin{split}
\mbox{\upshape I}&\leq\frac{C}{\sigma}\sum_i
\left(\ell(Q_i)^{\alpha}\int_{Q_i}|f(z)|\,dz\times\sum_{j=1}^\infty\frac{1}{[2^{j-1}\ell(Q_i)]^{n+\alpha}}\int_{2^{j+1}Q_i\backslash 2^{j}Q_i}\big|b(x)-b_{2^{j+1}Q_i}\big|\cdot w(x)\,dx\right).
\end{split}
\end{equation*}
Since $w\in A_1$, we know that there exists a number $r>1$ such that $w\in RH_r$. It then follows from H\"older's inequality, the John--Nirenberg's inequality(\cite{john}) and (\ref{general weights}) that
\begin{align}
\int_{2^{j+1}Q_i}\big|b(x)-b_{2^{j+1}Q_i}\big|\cdot w(x)\,dx&\leq\left(\int_{2^{j+1}Q_i}\big|b(x)-b_{2^{j+1}Q_i}\big|^{r'}dx\right)^{1/{r'}}\left(\int_{2^{j+1}Q_i}w(x)^rdx\right)^{1/r}\notag\\
&\leq C\|b\|_*\cdot w\big(2^{j+1}Q_i\big)\notag\\
&\leq C\|b\|_*\cdot(2^{j+1})^nw\big(Q_i\big).
\end{align}
Hence
\begin{equation*}
\begin{split}
\mbox{\upshape I}&\leq
\frac{C\cdot\|b\|_*}{\sigma}\sum_i\left(\int_{Q_i}|f(z)|\,dz\times\sum_{j=1}^\infty\frac{(2^{j+1})^nw\big(Q_i\big)}{(2^{j-1})^{n+\alpha}|Q_i|}\right)\\
&\leq\frac{C}{\sigma}\sum_i\left(\frac{w\big(Q_i\big)}{|Q_i|}\cdot\int_{Q_i}|f(z)|\,dz\times\sum_{j=1}^\infty\frac{1}{2^{j\alpha}}\right)\\
&\leq\frac{C}{\sigma}\sum_i\underset{z\in Q_i}{\mbox{ess\,inf}}\,w(z)\int_{Q_i}|f(z)|\,dz\\
&\leq \frac{C}{\sigma}\int_{\bigcup_i Q_i}|f(z)|w(z)\,dz\leq C\int_{\mathbb R^n}\frac{|f(z)|}{\sigma}\cdot w(z)\,dz\\
&\leq C\int_{\mathbb R^n}\Phi\left(\frac{|f(z)|}{\sigma}\right)\cdot w(z)\,dz.
\end{split}
\end{equation*}
For the term II, since $b\in BMO(\mathbb R^n)$, then a simple calculation shows that
\begin{equation}\label{BMO}
\big|b_{2^{j+1}Q_i}-b_{Q_i}\big|\leq C\cdot(j+1)\|b\|_*.
\end{equation}
This estimate (\ref{BMO}) together with the inequality (\ref{general weights}) implies that
\begin{equation*}
\begin{split}
\mbox{\upshape II}&\leq\frac{C\cdot\|b\|_*}{\sigma}\sum_i\left(\ell(Q_i)^{\alpha}\int_{Q_i}|f(z)|\,dz\times\sum_{j=1}^\infty
\big(j+1\big)\cdot\frac{w\big(2^{j+1}Q_i\big)}{[2^{j-1}\ell(Q_i)]^{n+\alpha}}\right)\\
&\leq\frac{C\cdot\|b\|_*}{\sigma}\sum_i\left(\int_{Q_i}|f(z)|\,dz\times\sum_{j=1}^\infty\big(j+1\big)\cdot\frac{(2^{j+1})^nw\big(Q_i\big)}{(2^{j-1})^{n+\alpha}|Q_i|}\right)\\
&\leq\frac{C}{\sigma}\sum_i\left(\frac{w\big(Q_i\big)}{|Q_i|}\cdot\int_{Q_i}|f(z)|\,dz\times\sum_{j=1}^\infty\frac{(j+1)}{2^{j\alpha}}\right)\\
&\leq\frac{C}{\sigma}\sum_i\left(\frac{w\big(Q_i\big)}{|Q_i|}\cdot\int_{Q_i}|f(z)|\,dz\right)\leq C\int_{\mathbb R^n}\Phi\left(\frac{|f(z)|}{\sigma}\right)\cdot w(z)\,dz.
\end{split}
\end{equation*}
On the other hand, by using the weighted weak-type (1,1) estimate of intrinsic square functions (see \cite{wilson2}), we have
\begin{equation*}
\begin{split}
I_6&\leq\frac{C}{\sigma}\int_{\mathbb R^n}\sum_i\big|b(x)-b_{Q_i}\big|\cdot|h_i(x)|\,w(x)\,dx\\
&=\frac{C}{\sigma}\sum_i\int_{Q_i}\big|b(x)-b_{Q_i}\big|\cdot|h_i(x)|\,w(x)\,dx\\
&\leq\frac{C}{\sigma}\sum_i\int_{Q_i}\big|b(x)-b_{Q_i}\big|\cdot|f(x)|\,w(x)\,dx\\
&+\frac{C}{\sigma}\sum_i\frac{1}{|Q_i|}\int_{Q_i}|f(y)|\,dy\times\int_{Q_i}\big|b(x)-b_{Q_i}\big|w(x)\,dx\\
&:=\mbox{\upshape III+IV}.
\end{split}
\end{equation*}
By the generalized H\"older's inequality with weight (\ref{weighted holder}), (\ref{weighted exp}) and (\ref{equiv norm with weight}), we can deduce that
\begin{equation*}
\begin{split}
\mbox{\upshape III}\leq&\frac{C}{\sigma}\sum_i w(Q_i)\cdot\frac{1}{w(Q_i)}\int_{Q_i}\big|b(x)-b_{Q_i}\big|\cdot|f(x)|\,w(x)\,dx\\
\leq&\frac{C}{\sigma}\sum_i w(Q_i)\cdot\big\|b-b_{Q_i}\big\|_{\exp L(w),Q_i}\big\|f\big\|_{L\log L(w),Q_i}\\
\leq&\frac{C\cdot\|b\|_*}{\sigma}\sum_i w(Q_i)\cdot\big\|f\big\|_{L\log L(w),Q_i}\\
=&\frac{C\cdot\|b\|_*}{\sigma}\sum_i w(Q_i)\cdot
\inf_{\eta>0}\left\{\eta+\frac{\eta}{w(Q_i)}\int_{Q_i}\Phi\left(\frac{|f(y)|}{\eta}\right)w(y)\,dy\right\}\\
\leq&\frac{C\cdot\|b\|_*}{\sigma}\sum_i w(Q_i)\cdot
\left\{\sigma+\frac{\sigma}{w(Q_i)}\int_{Q_i}\Phi\left(\frac{|f(y)|}{\sigma}\right)w(y)\,dy\right\}\\
\leq& C\left\{\sum_i w(Q_i)+\sum_i\int_{Q_i}\Phi\left(\frac{|f(y)|}{\sigma}\right)w(y)\,dy\right\}\\
\leq& C\int_{\mathbb R^n}\Phi\left(\frac{|f(y)|}{\sigma}\right)\cdot w(y)\,dy.
\end{split}
\end{equation*}
Arguing as in the proof of (3.8), we find that
\begin{equation*}
\begin{split}
\int_{Q_i}\big|b(x)-b_{Q_i}\big|w(x)\,dx&\leq\left(\int_{Q_i}\big|b(x)-b_{Q_i}\big|^{r'}dx\right)^{1/{r'}}\left(\int_{Q_i}w(x)^rdx\right)^{1/r}\\
&\leq C\|b\|_*\cdot w(Q_i).
\end{split}
\end{equation*}
Therefore
\begin{equation*}
\begin{split}
\mbox{\upshape IV}\leq&\frac{C}{\sigma}\sum_i\frac{w(Q_i)}{|Q_i|}\int_{Q_i}|f(y)|\,dy\leq\frac{C}{\sigma}\sum_i\underset{y\in Q_i}{\mbox{ess\,inf}}\,w(y)\int_{Q_i}|f(y)|\,dy\\
\leq&\frac{C}{\sigma}\int_{\bigcup_i Q_i}|f(y)|w(y)\,dy\leq C\int_{\mathbb R^n}\frac{|f(y)|}{\sigma}\cdot w(y)\,dy\\
\leq& C\int_{\mathbb R^n}\Phi\left(\frac{|f(y)|}{\sigma}\right)\cdot w(y)\,dy.
\end{split}
\end{equation*}
Summing up all the above estimates, we get the desired result.
\end{proof}

\section{Proof of Theorem \ref{mainthm:3}}

In order to prove the main theorem of this section, we will need the following estimates which were established by the author in \cite{wang5}.

\newtheorem{prop}[theorem]{Proposition}

\begin{prop}
Let $w\in A_1$ and $0<\alpha\le1$. Then for any $j\in\mathbb Z_+$, we have
\begin{equation*}
\big\|\mathcal S_{\alpha,2^j}(f)\big\|_{L^2_w}\leq C\cdot2^{jn/2}\big\|\mathcal S_\alpha(f)\big\|_{L^2_w}.
\end{equation*}
\end{prop}

\begin{prop}
Let $w\in A_1$, $0<\alpha\le1$ and $2<q<\infty$. Then for any $j\in\mathbb Z_+$, we have
\begin{equation*}
\big\|\mathcal S_{\alpha,2^j}(f)\big\|_{L^q_w}\leq C\cdot2^{jn/2}\big\|\mathcal S_\alpha(f)\big\|_{L^q_w}.
\end{equation*}
\end{prop}

\begin{prop}
Let $w\in A_1$, $0<\alpha\le1$ and $1<q<2$. Then for any $j\in\mathbb Z_+$, we have
\begin{equation*}
\big\|\mathcal S_{\alpha,2^j}(f)\big\|_{L^q_w}\leq C\cdot2^{jn/q}\big\|\mathcal S_\alpha(f)\big\|_{L^q_w}.
\end{equation*}
\end{prop}
Moreover, from the definition of $\mathcal G^*_{\lambda,\alpha}$($\lambda>3$), we readily see that
\begin{align}\label{g(lambda)}
\left|\mathcal G^*_{\lambda,\alpha}(f)(x)\right|^2=&\iint_{\mathbb R^{n+1}_+}\left(\frac{t}{t+|x-y|}\right)^{\lambda n}\Big(A_\alpha(f)(y,t)\Big)^2\frac{dydt}{t^{n+1}}\notag\\
=&\int_0^\infty\int_{|x-y|<t}\left(\frac{t}{t+|x-y|}\right)^{\lambda n}\Big(A_\alpha(f)(y,t)\Big)^2\frac{dydt}{t^{n+1}}\notag\\
&+\sum_{j=1}^\infty\int_0^\infty\int_{2^{j-1}t\le|x-y|<2^jt}\left(\frac{t}{t+|x-y|}\right)^{\lambda n}\Big(A_\alpha(f)(y,t)\Big)^2\frac{dydt}{t^{n+1}}\notag\\
\le&\, C\bigg[\mathcal S_\alpha(f)(x)^2+\sum_{j=1}^\infty 2^{-j\lambda n}\mathcal S_{\alpha,2^j}(f)(x)^2\bigg].
\end{align}
Thus, by applying Propositions 4.1--4.3, the $L^q_w$-boundedness of $\mathcal S_\alpha$(see \cite{wilson2}) and the above inequality (\ref{g(lambda)}), we obtain that for $1<q<\infty$ and $w\in A_1$,
\begin{align}\label{g(lambda) bounds}
\big\|\mathcal G^*_{\lambda,\alpha}(f)\big\|_{L^q_{w}}
&\le C\Bigg(\big\|\mathcal S_\alpha(f)\big\|_{L^q_{w}}+\sum_{j=1}^\infty 2^{-\frac{j\lambda n}{2}}\big\|\mathcal S_{\alpha,2^j}(f)\big\|_{L^q_{w}}\Bigg)\notag\\
&\le C\Bigg(\big\|\mathcal S_\alpha(f)\big\|_{L^q_{w}}+\sum_{j=1}^\infty 2^{-\frac{j\lambda n}{2}}\cdot\big[2^{\frac{jn}{2}}+2^{\frac{jn}{q}}\big]\big\|\mathcal S_{\alpha}(f)\big\|_{L^q_{w}}\Bigg)\notag\\
&\le C\big\|f\big\|_{L^q_{w}}\Bigg(1+\sum_{j=1}^\infty2^{-\frac{j\lambda n}{2}}
\cdot\big[2^{\frac{jn}{2}}+2^{\frac{jn}{q}}\big]\Bigg)\notag\\
&\le C\big\|f\big\|_{L^q_{w}},
\end{align}
where the last inequality holds under the assumption $\lambda>3>\max\{1,2/q\}$ when $1<q<\infty$. In addition, for a given real-valued function $b\in BMO(\mathbb R^n)$, as before, we can also prove that
\begin{equation}\label{4.3}
\big|\big[b,\mathcal G^*_{\lambda,\alpha}\big](f)(x)\big|\le\frac{1}{2\pi}\int_0^{2\pi}
\mathcal G^*_{\lambda,\alpha}\big(e^{-e^{i\theta}b}\cdot f\big)(x)\cdot e^{\cos\theta\cdot b(x)}\,d\theta.
\end{equation}
Taking into account the inequalities (\ref{g(lambda) bounds}) and (\ref{4.3}), and following along the same arguments used in \cite{ding}, we can also show the following:

\begin{theorem}\label{commutator thm2}
Let $0<\alpha\le1$, $1<q<\infty$ and $w\in A_1$. Suppose that $\lambda>3$, then the commutator $\big[b,\mathcal G^*_{\lambda,\alpha}\big]$ are bounded from $L^q_w(\mathbb R^n)$ into itself whenever $b\in BMO(\mathbb R^n)$.
\end{theorem}

In \cite{wang2}, we have established the weighted weak-type (1,1) estimate of $\mathcal G^*_{\lambda,\alpha}$ on $L^1_w(\mathbb R^n)$. More specifically, we obtained
\begin{theorem}\label{g(lambda)weak estimate}
Let $0<\alpha\le1$, $w\in A_1$ and $\lambda>{(3n+2\alpha)}/n$. Then for any given $\sigma>0$, there exists a constant $C>0$ independent of $f$ and $\sigma$ such that
\begin{equation*}
w\Big(\Big\{x\in\mathbb R^n:\big|\mathcal G^*_{\lambda,\alpha}(f)(x)\big|>\sigma\Big\}\Big)\le\frac{C}{\sigma}\int_{\mathbb R^n}|f(x)|w(x)\,dx.
\end{equation*}
\end{theorem}

\begin{proof}[Proof of Theorem $\ref{mainthm:3}$]
For any fixed $\sigma>0$, as before, we again perform the Calder\'on--Zygmund decomposition of $f$ at the level $\sigma$ to obtain a sequence of disjoint non-overlapping dyadic cubes $\{Q_i\}$ such that the following property holds (see \cite{stein})
\begin{equation}\label{decomposition2}
\sigma<\frac{1}{|Q_i|}\int_{Q_i}|f(y)|\,dy< 2^n\sigma.
\end{equation}
Setting $E=\bigcup_i Q_i$. Now we decompose $f(x)=g(x)+h(x)$, where

$g(x)=f(x)$ when $x\in E^c$, and $g(x)=\frac{1}{|Q_i|}\int_{Q_i}|f(y)|\,dy$ when $x\in Q_i$. Then
\begin{equation*}
h(x)=f(x)-g(x)=\sum_i h_i(x),
\end{equation*}
with $h_i(x)=h(x)\chi_{Q_i}(x)$. Clearly, by the above decomposition, we get supp\,$h_i\subseteq Q_i$, $\int_{Q_i}h_i(x)\,dx=0$ and $\|h_i\|_{L^1}\le 2\int_{Q_i}|f(x)|\,dx$ . Note that $\big|\big[b,\mathcal G^*_{\lambda,\alpha}\big](f)(x)\big|\leq\big|\big[b,\mathcal G^*_{\lambda,\alpha}\big](g)(x)\big|+\big|\big[b,\mathcal G^*_{\lambda,\alpha}\big](h)(x)\big|$, then we have
\begin{equation*}
\begin{split}
&w\big(\big\{x\in\mathbb R^n:\big|[b,\mathcal G^*_{\lambda,\alpha}](f)(x)\big|>\sigma\big\}\big)\\
\leq&
w\big(\big\{x\in\mathbb R^n:\big|[b,\mathcal G^*_{\lambda,\alpha}](g)(x)\big|>\sigma/2\big\}\big)
+w\big(\big\{x\in\mathbb R^n:\big|[b,\mathcal G^*_{\lambda,\alpha}](h)(x)\big|>\sigma/2\big\}\big)\\
:=&J_1+J_2.
\end{split}
\end{equation*}
Let us start with the term $J_1$. By using Chebyshev's inequality, Theorem \ref{commutator thm2} and the inequality (\ref{g}), we obtain
\begin{equation*}
\begin{split}
J_1&\leq \frac{4}{\sigma^2}\cdot\Big\|\big[b,\mathcal G^*_{\lambda,\alpha}\big](g)\Big\|^2_{L^2_w}\leq \frac{C}{\sigma^2}\cdot\big\|g\big\|^2_{L^2_w}\\
&\leq\frac{C}{\sigma^2}\cdot\sigma\int_{\mathbb R^n}|f(x)|w(x)\,dx\\
&\leq C\int_{\mathbb R^n}\Phi\left(\frac{|f(x)|}{\sigma}\right)\cdot w(x)\,dx.
\end{split}
\end{equation*}
To estimate the other term $J_2$, as before, we also let $Q_i^*=2\sqrt n Q_i$ be the cube concentric with $Q_i$ such that $\ell(Q_i^*)=(2\sqrt n)\ell(Q_i)$. Then we can further split $J_2$ into two parts as follows.
\begin{equation*}
\begin{split}
J_2\le&\,w\Big(\Big\{x\in \bigcup_i Q_i^*:\Big|\big[b,\mathcal G^*_{\lambda,\alpha}\big](h)(x)\Big|>\sigma/2\Big\}\Big)\\
&+w\Big(\Big\{x\notin \bigcup_i Q_i^*:\Big|\big[b,\mathcal G^*_{\lambda,\alpha}\big](h)(x)\Big|>\sigma/2\Big\}\Big)\\
:=&\,J_3+J_4.
\end{split}
\end{equation*}
The part of the argument involving $J_3$ proceeds as in Theorem \ref{mainthm:1},
\begin{equation*}
J_3\leq\sum_i w\big(Q_i^*\big)\le C\sum_i w(Q_i)\leq C\int_{\mathbb R^n}\Phi\left(\frac{|f(x)|}{\sigma}\right)\cdot w(x)\,dx.
\end{equation*}
By the previous estimate (\ref{commutator estimate}), we thus obtain
\begin{equation*}
\begin{split}
\big|\big[b,\mathcal G^*_{\lambda,\alpha}\big](h)(x)\big|&\leq\sum_i\big|b(x)-b_{Q_i}\big|\cdot\mathcal G^*_{\lambda,\alpha}(h_i)(x)\\
&+\left(\iint_{\Gamma(x)}\sup_{\varphi\in{\mathcal C}_\alpha}\bigg|\int_{\mathbb R^n}\big[b(z)-b_{Q_i}\big]\varphi_t(y-z)\cdot\sum_i h_i(z)\,dz\bigg|^2\frac{dydt}{t^{n+1}}\right)^{1/2}\\
&=\sum_i\big|b(x)-b_{Q_i}\big|\cdot\mathcal G^*_{\lambda,\alpha}(h_i)(x)+\mathcal G^*_{\lambda,\alpha}\bigg(\sum_i[b-b_{Q_i}]h_i\bigg)(x).
\end{split}
\end{equation*}
Therefore
\begin{equation*}
\begin{split}
J_4\leq & w\bigg(\bigg\{x\notin \bigcup_i Q_i^*:\sum_i\big|b(x)-b_{Q_i}\big|\cdot\mathcal G^*_{\lambda,\alpha}(h_i)(x)>\sigma/4\bigg\}\bigg)\\
&+ w\bigg(\bigg\{x\notin \bigcup_i Q_i^*:\mathcal G^*_{\lambda,\alpha}\bigg(\sum_i[b-b_{Q_i}]h_i\bigg)(x)>\sigma/4\bigg\}\bigg)\\
:=&J_5+J_6.
\end{split}
\end{equation*}
It follows directly from the Chebyshev's inequality that
\begin{equation*}
\begin{split}
J_5&\leq\frac{\,4\,}{\sigma}\int_{\mathbb R^n\backslash\bigcup_i Q_i^*}\bigg|\sum_i\big|b(x)-b_{Q_i}\big|\cdot\mathcal G^*_{\lambda,\alpha}(h_i)(x)\bigg|w(x)\,dx\\
&\leq\frac{\,4\,}{\sigma}\sum_i\left(\int_{(Q_i^*)^c}\big|b(x)-b_{Q_i}\big|\cdot\mathcal G^*_{\lambda,\alpha}(h_i)(x)\,w(x)\,dx\right).\\
\end{split}
\end{equation*}
We also denote the center of $Q_i$ by $c_i$. In the proof of Theorem \ref{mainthm:1}, we have already shown that
\begin{equation}
\big|\mathcal S_{\alpha}(h_i)(x)\big|\leq C\cdot\frac{\ell(Q_i)^{\alpha}}{|x-c_i|^{n+\alpha}}\bigg(\int_{Q_i}|f(z)|\,dz\bigg).
\end{equation}
Below we will give the pointwise estimates of $\big|\mathcal S_{\alpha,2^j}(h_i)(x)\big|$ for $j=1,2,\ldots$. Notice that for any $z\in Q_i$ and $x\in (Q^*_i)^c$, we get $|z-c_i|<\frac{|x-c_i|}{2}$. Thus, for all $(y,t)\in\Gamma_{2^j}(x)$ and $|z-y|\le t$ with $z\in Q_i$, we can deduce that
\begin{equation}\label{2(j)t}
t+2^jt\ge|x-y|+|y-z|\ge|x-z|\ge|x-c_i|-|z-c_i|\ge\frac{|x-c_i|}{2}.
\end{equation}
Hence, for any $x\in (Q^*_i)^c$, by the inequalities (\ref{kernel estimate1}) and (\ref{2(j)t}), we obtain that for $j=1,2,\ldots$,
\begin{equation*}
\begin{split}
\big|\mathcal S_{\alpha,2^j}(h_i)(x)\big|&=\left(\iint_{\Gamma_{2^j}(x)}\bigg(\sup_{\varphi\in{\mathcal C}_\alpha}\big|(\varphi_t*{h_i})(y)\big|\bigg)^2\frac{dydt}{t^{n+1}}\right)^{1/2}\\
&\leq C\cdot\ell(Q_i)^{\alpha}\bigg(\int_{Q_i}|h_i(z)|\,dz\bigg)\left(\int_{\frac{|x-c_i|}{2^{j+2}}}^\infty
\int_{|y-x|<2^jt}\frac{dydt}{t^{2(n+\alpha)+n+1}}\right)^{1/2}\\
&\leq C\cdot2^{{jn}/2}\ell(Q_i)^{\alpha}\bigg(\int_{Q_i}|h_i(z)|\,dz\bigg)
\left(\int_{\frac{|x-c_i|}{2^{j+2}}}^\infty\frac{dt}{t^{2(n+\alpha)+1}}\right)^{1/2}\\
&\leq C\cdot2^{j(3n+2\alpha)/2}\frac{\ell(Q_i)^{\alpha}}{|x-c_i|^{n+\alpha}}\bigg(\int_{Q_i}|f(z)|\,dz\bigg).
\end{split}
\end{equation*}
Therefore, by using the pointwise estimate we just derived above and the inequality (\ref{g(lambda)}),
\begin{equation*}
\begin{split}
\left|\mathcal G^*_{\lambda,\alpha}(h_i)(x)\right|&\leq C\bigg[\big|\mathcal S_\alpha(h_i)(x)\big|
+\sum_{j=1}^\infty 2^{{-j\lambda n}/2}\big|\mathcal S_{\alpha,2^j}(h_i)(x)\big|\bigg]\\
&\leq C\cdot\frac{\ell(Q_i)^{\alpha}}{|x-c_i|^{n+\alpha}}\bigg(\int_{Q_i}|f(z)|\,dz\bigg)\times\left(1+\sum_{j=1}^\infty 2^{{-j\lambda n}/2}\cdot2^{j(3n+2\alpha)/2}\right)\\
&\leq C\cdot\frac{\ell(Q_i)^{\alpha}}{|x-c_i|^{n+\alpha}}\bigg(\int_{Q_i}|f(z)|\,dz\bigg),
\end{split}
\end{equation*}
where the last inequality is due to our assumption $\lambda>{(3n+2\alpha)}/n$. Consequently,
\begin{equation*}
J_5\leq\frac{C}{\sigma}\sum_i
\left(\ell(Q_i)^{\alpha}\int_{Q_i}|f(z)|\,dz\times\int_{(Q_i^*)^c}\big|b(x)-b_{Q_i}\big|\cdot\frac{w(x)}{|x-c_i|^{n+\alpha}}dx\right).
\end{equation*}
Following along the same lines as in Theorem \ref{mainthm:1}, we can also show
\begin{equation*}
J_5\leq C\int_{\mathbb R^n}\Phi\left(\frac{|f(x)|}{\sigma}\right)\cdot w(x)\,dx.
\end{equation*}
On the other hand, by using the weighted weak-type (1,1) estimate of $\mathcal G^*_{\lambda,\alpha}$(see Theorem \ref{g(lambda)weak estimate}), we have
\begin{equation*}
J_6\leq\frac{C}{\sigma}\int_{\mathbb R^n}\sum_i\big|b(x)-b_{Q_i}\big|\cdot|h_i(x)|\,w(x)\,dx.
\end{equation*}
The rest of the proof is exactly the same as that of Theorem \ref{mainthm:1}, and we finally obtain
\begin{equation*}
J_6\leq C\int_{\mathbb R^n}\Phi\left(\frac{|f(x)|}{\sigma}\right)\cdot w(x)\,dx.
\end{equation*}
Collecting all these estimates, we get the desired estimate.
\end{proof}

\section{Proofs of Theorems \ref{mainthm:4}, \ref{mainthm:5} and \ref{mainthm:6}}

\begin{proof}[Proofs of Theorems $\ref{mainthm:4}$ and $\ref{mainthm:5}$]
We will only give the proof of Theorem $\ref{mainthm:4}$ here, because the proof of Theorem $\ref{mainthm:5}$ is essentially the same. Fix a ball $B=B(x_0,r_B)\subseteq\mathbb R^n$ and decompose $f=f_1+f_2$, where $f_1=f\cdot\chi_{_{2B}}$, $\chi_{_{2B}}$ denotes the characteristic function of $2B=B(x_0,2r_B)$. For any $0<\kappa<1$, $w\in A_1$ and any given $\sigma>0$, one writes
\begin{equation*}
\begin{split}
&\frac{1}{w(B)^\kappa}\cdot w\big(\big\{x\in B:\big|[b,\mathcal S_\alpha](f)(x)\big|>\sigma\big\}\big)\\
\leq &\frac{1}{w(B)^\kappa}\cdot w\big(\big\{x\in B:\big|[b,\mathcal S_\alpha](f_1)(x)\big|>\sigma/2\big\}\big)
+\frac{1}{w(B)^\kappa}\cdot w\big(\big\{x\in B:\big|[b,\mathcal S_\alpha](f_2)(x)\big|>\sigma/2\big\}\big)\\
:=&I_1+I_2.
\end{split}
\end{equation*}
Using Theorem \ref{mainthm:1} and the inequality (\ref{weights}), we get
\begin{equation*}
\begin{split}
I_1&\leq C\cdot\frac{1}{w(B)^\kappa}\int_{\mathbb R^n}\Phi\left(\frac{|f_1(x)|}{\sigma}\right)\cdot w(x)\,dx\\
&= C\cdot\frac{1}{w(B)^\kappa}
\int_{2B}\Phi\left(\frac{|f(x)|}{\sigma}\right)\cdot w(x)\,dx\\
&= C\cdot\frac{w(2B)^\kappa}{w(B)^\kappa}\cdot\frac{1}{w(2B)^\kappa}
\int_{2B}\Phi\left(\frac{|f(x)|}{\sigma}\right)\cdot w(x)\,dx\\
&\leq C\cdot\sup_B\left\{\frac{1}{w(B)^\kappa}
\int_{B}\Phi\left(\frac{|f(x)|}{\sigma}\right)\cdot w(x)\,dx\right\}.
\end{split}
\end{equation*}
For any $x\in B$, we can easily check that
\begin{equation*}
\big|\big[b,\mathcal S_\alpha\big](f_2)(x)\big|\leq \big|b(x)-b_{B}\big|\cdot\mathcal S_\alpha(f_2)(x)+\mathcal S_\alpha\Big([b-b_{B}]f_2\Big)(x).
\end{equation*}
So we have
\begin{equation*}
\begin{split}
I_2\leq&\frac{1}{w(B)^\kappa}\cdot w\big(\big\{x\in B:\big|b(x)-b_{B}\big|\cdot\mathcal S_\alpha(f_2)(x)>\sigma/4\big\}\big)\\
&+\frac{1}{w(B)^\kappa}\cdot w\big(\big\{x\in B:\big|\mathcal S_\alpha\Big([b-b_{B}]f_2\Big)(x)\big|>\sigma/4\big\}\big)\\
:=&I_3+I_4.
\end{split}
\end{equation*}
For the term $I_3$, for all $0<\alpha\leq1$ and $x\in B$, it was proved by the author \cite{wang1} that
\begin{equation}\label{S(f2)}
\big|\mathcal S_\alpha(f_2)(x)\big|\leq C\sum_{j=1}^\infty\frac{1}{|2^{j+1}B|}\int_{2^{j+1}B}|f(z)|\,dz.
\end{equation}
Since $w\in A_1$, then there exists a number $r>1$ such that $w\in RH_r$. Hence, by using the above pointwise estimate (\ref{S(f2)}), Chebyshev's inequality together with H\"older's inequality and John--Nirenberg's inequality (see \cite{john}), we conclude that
\begin{equation*}
\begin{split}
I_3&\leq\frac{1}{w(B)^\kappa}\cdot\frac{\,4\,}{\sigma}\int_B\big|b(x)-b_{B}\big|\cdot\mathcal S_\alpha(f_2)(x)w(x)\,dx\\
&\leq C\sum_{j=1}^\infty\frac{1}{|2^{j+1}B|}\int_{2^{j+1}B}\frac{|f(z)|}{\sigma}\,dz\\
&\times\frac{1}{w(B)^\kappa}\cdot\left(\int_{B}\big|b(x)-b_{B}\big|^{r'}dx\right)^{1/{r'}}\left(\int_{B}w(x)^rdx\right)^{1/r}\\
&\leq C\sum_{j=1}^\infty\frac{1}{|2^{j+1}B|}\int_{2^{j+1}B}\frac{|f(z)|}{\sigma}\,dz\times w(B)^{1-\kappa}\\
&\leq C\sum_{j=1}^\infty\frac{1}{w(2^{j+1}B)}\int_{2^{j+1}B}\frac{|f(z)|}{\sigma}\cdot w(z)\,dz\times w(B)^{1-\kappa}\\
&\leq C\cdot\sup_B\left\{\frac{1}{w(B)^\kappa}
\int_{B}\Phi\left(\frac{|f(z)|}{\sigma}\right)\cdot w(z)\,dz\right\}\times\sum_{j=1}^\infty\frac{w(B)^{1-\kappa}}{w(2^{j+1}B)^{1-\kappa}}.
\end{split}
\end{equation*}
Since $w\in A_1\subset A_\infty$, by the inequality (\ref{compare}), we get
\begin{align}\label{<C}
\sum_{j=1}^\infty\frac{w(B)^{1-\kappa}}{w(2^{j+1}B)^{1-\kappa}}&\leq C\sum_{j=1}^\infty\left(\frac{|B|}{|2^{j+1}B|}\right)^{\delta(1-\kappa)}\notag\\
&\leq C\sum_{j=1}^\infty\left(\frac{1}{2^{jn}}\right)^{\delta(1-\kappa)}\leq C,
\end{align}
which in turn gives that
\begin{equation*}
I_3\leq C\cdot\sup_B\left\{\frac{1}{w(B)^\kappa}
\int_{B}\Phi\left(\frac{|f(z)|}{\sigma}\right)\cdot w(z)\,dz\right\}.
\end{equation*}
Similar to the proof of (\ref{S(f2)}), for all $0<\alpha\leq1$ and all $x\in B$, we can show the following pointwise estimate as well.
\begin{equation}\label{[b,S](f2)}
\Big|\mathcal S_\alpha\Big([b-b_{B}]f_2\Big)(x)\Big|\leq C\sum_{j=1}^\infty\frac{1}{|2^{j+1}B|}\int_{2^{j+1}B}\big|b(z)-b_{B}\big|\cdot\big|f(z)\big|\,dz.
\end{equation}
Applying the above pointwise estimate (\ref{[b,S](f2)}) and Chebyshev's inequality, we have
\begin{equation*}
\begin{split}
I_4&\leq\frac{1}{w(B)^\kappa}\cdot\frac{\,4\,}{\sigma}\int_B\Big|\mathcal S_\alpha\Big([b-b_{B}]f_2\Big)(x)\Big|w(x)\,dx\\
&\leq\frac{w(B)}{w(B)^\kappa}\cdot\frac{C}{\sigma}\sum_{j=1}^\infty\frac{1}{|2^{j+1}B|}\int_{2^{j+1}B}
\big|b(z)-b_{B}\big|\cdot\big|f(z)\big|\,dz\\
\end{split}
\end{equation*}
\begin{equation*}
\begin{split}
&\leq\frac{w(B)}{w(B)^\kappa}\cdot\frac{C}{\sigma}\sum_{j=1}^\infty\frac{1}{|2^{j+1}B|}\int_{2^{j+1}B}
\big|b(z)-b_{2^{j+1}B}\big|\cdot\big|f(z)\big|\,dz\\
&+\frac{w(B)}{w(B)^\kappa}\cdot\frac{C}{\sigma}\sum_{j=1}^\infty\frac{1}{|2^{j+1}B|}\int_{2^{j+1}B}
\big|b_{2^{j+1}B}-b_B\big|\cdot\big|f(z)\big|\,dz\\
&:=I_5+I_6.
\end{split}
\end{equation*}
For the term $I_5$, observe that for any $a,b>0$, $\Phi(a\cdot b)\leq\Phi(a)\cdot\Phi(b)$ when $\Phi(t)=t(1+\log^+t)$. We then use the generalized H\"older's inequality with weight (\ref{weighted holder}), (\ref{weighted exp}) and (\ref{equiv norm with weight}) together with (\ref{<C}) to obtain
\begin{equation*}
\begin{split}
I_5&\leq \frac{C}{\sigma}\cdot w(B)^{1-\kappa}\sum_{j=1}^\infty\frac{1}{w(2^{j+1}B)}\int_{2^{j+1}B}
\big|b(z)-b_{2^{j+1}B}\big|\cdot\big|f(z)\big|w(z)\,dz\\
&\leq \frac{C}{\sigma}\cdot w(B)^{1-\kappa}\sum_{j=1}^\infty\big\|b-b_{2^{j+1}B}\big\|_{\exp L(w),2^{j+1}B}\big\|f\big\|_{L\log L(w),2^{j+1}B}\\
&\leq \frac{C\|b\|_*}{\sigma}\cdot w(B)^{1-\kappa}\sum_{j=1}^\infty
\inf_{\eta>0}\left\{\eta+\frac{\eta}{w(2^{j+1}B)}\int_{2^{j+1}B}\Phi\left(\frac{|f(z)|}{\eta}\right)w(z)\,dz\right\}\\
&\leq \frac{C\|b\|_*}{\sigma}\cdot w(B)^{1-\kappa}\sum_{j=1}^\infty
\left\{\frac{\sigma}{w(2^{j+1}B)^{1-\kappa}}+\frac{\sigma}{w(2^{j+1}B)}\int_{2^{j+1}B}\Phi\left(\frac{|f(z)|}{\sigma}\right)w(z)\,dz\right\}\\
&\leq C\|b\|_*\cdot\left[1+\sup_B\left\{\frac{1}{w(B)^\kappa}
\int_{B}\Phi\left(\frac{|f(z)|}{\sigma}\right)\cdot w(z)\,dz\right\}\right]\times\sum_{j=1}^\infty\frac{w(B)^{1-\kappa}}{w(2^{j+1}B)^{1-\kappa}}\\
&\leq C\cdot\sup_B\left\{\frac{1}{w(B)^\kappa}
\int_{B}\Phi\left(\frac{|f(z)|}{\sigma}\right)\cdot w(z)\,dz\right\}.
\end{split}
\end{equation*}
For the last term $I_6$ we proceed as follows. By the inequality (\ref{BMO}) again, we get
\begin{equation*}
\begin{split}
I_6&\leq C\cdot w(B)^{1-\kappa}\sum_{j=1}^\infty(j+1)\|b\|_*\cdot\frac{1}{|2^{j+1}B|}\int_{2^{j+1}B}\frac{|f(z)|}{\sigma}\,dz\\
&\leq C\cdot w(B)^{1-\kappa}\sum_{j=1}^\infty(j+1)\|b\|_*\cdot\frac{1}{w(2^{j+1}B)}\int_{2^{j+1}B}\frac{|f(z)|}{\sigma}\cdot w(z)\,dz\\
&\leq C\cdot\sup_B\left\{\frac{1}{w(B)^\kappa}
\int_{B}\Phi\left(\frac{|f(z)|}{\sigma}\right)\cdot w(z)\,dz\right\}\times\sum_{j=1}^\infty(j+1)\cdot\frac{w(B)^{1-\kappa}}{w(2^{j+1}B)^{1-\kappa}}.
\end{split}
\end{equation*}
Since $w\in A_1\subset A_\infty$, by using the inequality (\ref{compare}) again, we have
\begin{align}
\sum_{j=1}^\infty(j+1)\cdot\frac{w(B)^{1-\kappa}}{w(2^{j+1}B)^{1-\kappa}}&\leq C\sum_{j=1}^\infty(j+1)\cdot\left(\frac{|B|}{|2^{j+1}B|}\right)^{\delta(1-\kappa)}\notag\\
&\leq C\sum_{j=1}^\infty(j+1)\cdot\left(\frac{1}{2^{(j+1)n}}\right)^{\delta(1-\kappa)}\leq C,
\end{align}
which implies
\begin{equation*}
I_6\leq C\cdot\sup_B\left\{\frac{1}{w(B)^\kappa}
\int_{B}\Phi\left(\frac{|f(z)|}{\sigma}\right)\cdot w(z)\,dz\right\}.
\end{equation*}
Summarizing the above discussions, we obtain the conclusion of the theorem.
\end{proof}

\begin{proof}[Proof of Theorem $\ref{mainthm:6}$]
Fix a ball $B=B(x_0,r_B)\subseteq\mathbb R^n$ and decompose $f=f_1+f_2$, where $f_1=f\cdot\chi_{_{2B}}$. For any $0<\kappa<1$, $w\in A_1$ and any given $\sigma>0$, we then write
\begin{equation*}
\begin{split}
&\frac{1}{w(B)^\kappa}\cdot w\big(\big\{x\in B:\big|[b,\mathcal G^*_{\lambda,\alpha}](f)(x)\big|>\sigma\big\}\big)\\
\leq &\frac{1}{w(B)^\kappa}\cdot w\big(\big\{x\in B:\big|[b,\mathcal G^*_{\lambda,\alpha}](f_1)(x)\big|>\sigma/2\big\}\big)
+\frac{1}{w(B)^\kappa}\cdot w\big(\big\{x\in B:\big|[b,\mathcal G^*_{\lambda,\alpha}](f_2)(x)\big|>\sigma/2\big\}\big)\\
:=&J_1+J_2.
\end{split}
\end{equation*}
Theorem \ref{mainthm:3} and the inequality (\ref{weights}) imply that
\begin{equation*}
\begin{split}
J_1&\leq C\cdot\frac{1}{w(B)^\kappa}\int_{\mathbb R^n}\Phi\left(\frac{|f_1(x)|}{\sigma}\right)\cdot w(x)\,dx\\
&= C\cdot\frac{1}{w(B)^\kappa}
\int_{2B}\Phi\left(\frac{|f(x)|}{\sigma}\right)\cdot w(x)\,dx\\
&= C\cdot\frac{w(2B)^\kappa}{w(B)^\kappa}\cdot\frac{1}{w(2B)^\kappa}
\int_{2B}\Phi\left(\frac{|f(x)|}{\sigma}\right)\cdot w(x)\,dx\\
&\leq C\cdot\sup_B\left\{\frac{1}{w(B)^\kappa}
\int_{B}\Phi\left(\frac{|f(x)|}{\sigma}\right)\cdot w(x)\,dx\right\}.
\end{split}
\end{equation*}
For any $x\in B$, we are able to verify that
\begin{equation*}
\big|\big[b,\mathcal G^*_{\lambda,\alpha}\big](f_2)(x)\big|\leq \big|b(x)-b_{B}\big|\cdot\mathcal G^*_{\lambda,\alpha}(f_2)(x)
+\mathcal G^*_{\lambda,\alpha}\Big([b-b_{B}]f_2\Big)(x).
\end{equation*}
So we have
\begin{equation*}
\begin{split}
J_2\leq&\frac{1}{w(B)^\kappa}\cdot w\big(\big\{x\in B:\big|b(x)-b_{B}\big|\cdot\mathcal G^*_{\lambda,\alpha}(f_2)(x)>\sigma/4\big\}\big)\\
&+\frac{1}{w(B)^\kappa}\cdot w\big(\big\{x\in B:\big|\mathcal G^*_{\lambda,\alpha}\Big([b-b_{B}]f_2\Big)(x)\big|>\sigma/4\big\}\big)\\
:=&J_3+J_4.
\end{split}
\end{equation*}
For the term $J_3$, for all $0<\alpha\leq1$, $x\in B$ and $j\in\mathbb Z_+$, it was also shown by the author \cite{wang1} that
\begin{equation}\label{Sj(f2)}
\big|\mathcal S_{\alpha,2^j}(f_2)(x)\big|\leq C\cdot 2^{{3jn}/2}\sum_{j=1}^\infty\frac{1}{|2^{j+1}B|}\int_{2^{j+1}B}|f(z)|\,dz.
\end{equation}
Hence, it follows from the inequalities (\ref{Sj(f2)}), (\ref{S(f2)}) and (\ref{g(lambda)}) that
\begin{align}\label{g(lambda)(f2)}
\left|\mathcal G^*_{\lambda,\alpha}(f_2)(x)\right|&\leq C\bigg[\big|\mathcal S_\alpha(f_2)(x)\big|
+\sum_{j=1}^\infty 2^{{-j\lambda n}/2}\big|\mathcal S_{\alpha,2^j}(f_2)(x)\big|\bigg]\notag\\
&\leq C\cdot\sum_{j=1}^\infty\frac{1}{|2^{j+1}B|}\int_{2^{j+1}B}|f(z)|\,dz\times\left(1+\sum_{j=1}^\infty 2^{{-j\lambda n}/2}\cdot2^{3jn/2}\right)\notag\\
&\leq C\cdot\sum_{j=1}^\infty\frac{1}{|2^{j+1}B|}\int_{2^{j+1}B}|f(z)|\,dz,
\end{align}
where the last inequality is due to our assumption $\lambda>{(3n+2\alpha)}/n>3$.
Hence, we can continue the estimate of $J_3$ in the same way as in Theorem \ref{mainthm:4}, and obtain
\begin{equation*}
\begin{split}
J_3&\leq\frac{1}{w(B)^\kappa}\cdot\frac{\,4\,}{\sigma}\int_B\big|b(x)-b_{B}\big|\cdot\mathcal G^*_{\lambda,\alpha}(f_2)(x)w(x)\,dx\\
&\leq C\sum_{j=1}^\infty\frac{1}{|2^{j+1}B|}\int_{2^{j+1}B}\frac{|f(z)|}{\sigma}\,dz
\times\frac{1}{w(B)^\kappa}\cdot\int_B\big|b(x)-b_{B}\big|w(x)\,dx\\
&\leq C\sum_{j=1}^\infty\frac{1}{w(2^{j+1}B)}\int_{2^{j+1}B}\frac{|f(z)|}{\sigma}\cdot w(z)\,dz\times w(B)^{1-\kappa}\\
&\leq C\cdot\sup_B\left\{\frac{1}{w(B)^\kappa}
\int_{B}\Phi\left(\frac{|f(z)|}{\sigma}\right)\cdot w(z)\,dz\right\}\times\sum_{j=1}^\infty\frac{w(B)^{1-\kappa}}{w(2^{j+1}B)^{1-\kappa}}\\
&\leq C\cdot\sup_B\left\{\frac{1}{w(B)^\kappa}
\int_{B}\Phi\left(\frac{|f(z)|}{\sigma}\right)\cdot w(z)\,dz\right\}.
\end{split}
\end{equation*}
For the term $J_4$, similar to the proof of (\ref{g(lambda)(f2)}), for all $0<\alpha\leq1$, all $x\in B$ and $\lambda>3$, we can show the following pointwise estimate as well.
\begin{equation}\label{[b,g](f2)}
\Big|\mathcal G^*_{\lambda,\alpha}\Big([b-b_{B}]f_2\Big)(x)\Big|\leq C\sum_{j=1}^\infty\frac{1}{|2^{j+1}B|}\int_{2^{j+1}B}\big|b(z)-b_{B}\big|\cdot\big|f(z)\big|\,dz.
\end{equation}
Following the same arguments as in the proof of Theorem \ref{mainthm:4} and using the pointwise estimate (\ref{[b,g](f2)}) and Chebyshev's inequality, we have eventually obtained
\begin{equation*}
\begin{split}
J_4&\leq\frac{1}{w(B)^\kappa}\cdot\frac{\,4\,}{\sigma}\int_B\Big|\mathcal G^*_{\lambda,\alpha}\Big([b-b_{B}]f_2\Big)(x)\Big|w(x)\,dx\\
&\leq\frac{w(B)}{w(B)^\kappa}\cdot\frac{C}{\sigma}\sum_{j=1}^\infty\frac{1}{|2^{j+1}B|}\int_{2^{j+1}B}
\big|b(z)-b_{B}\big|\cdot\big|f(z)\big|\,dz\\
&\leq C\cdot\sup_B\left\{\frac{1}{w(B)^\kappa}
\int_{B}\Phi\left(\frac{|f(z)|}{\sigma}\right)\cdot w(z)\,dz\right\}.
\end{split}
\end{equation*}
Combining all the above estimates, we are done.
\end{proof}

\section{Proofs of Theorems \ref{mainthm:7}, \ref{mainthm:8} and \ref{mainthm:9}}

\begin{proof}[Proofs of Theorems $\ref{mainthm:7}$ and $\ref{mainthm:8}$]
Again we will only give the proof of Theorem $\ref{mainthm:7}$ here. Theorem $\ref{mainthm:8}$ can be dealt with similarly. For any ball $B=B(x_0,r)\subseteq\mathbb R^n$ with $x_0\in\mathbb R^n$ and $r>0$, we set $f=f_1+f_2$, where $f_1=f\cdot\chi_{_{2B}}$. Then for each fixed $\sigma>0$, we have
\begin{equation*}
\begin{split}
&\frac{1}{\Theta(r)}\cdot\big|\big\{x\in B:\big|[b,\mathcal S_\alpha](f)(x)\big|>\sigma\big\}\big|\\
\leq &\frac{1}{\Theta(r)}\cdot \big|\big\{x\in B:\big|[b,\mathcal S_\alpha](f_1)(x)\big|>\sigma/2\big\}\big|
+\frac{1}{\Theta(r)}\cdot \big|\big\{x\in B:\big|[b,\mathcal S_\alpha](f_2)(x)\big|>\sigma/2\big\}\big|\\
:=&K_1+K_2.
\end{split}
\end{equation*}
We consider the term $K_1$ first. Theorem \ref{maincor:1} and the inequality (\ref{doubling}) imply that
\begin{equation*}
\begin{split}
K_1&\leq C\cdot\frac{1}{\Theta(r)}\int_{\mathbb R^n}\Phi\left(\frac{|f_1(x)|}{\sigma}\right)dx\\
&= C\cdot\frac{1}{\Theta(r)}
\int_{2B}\Phi\left(\frac{|f(x)|}{\sigma}\right)dx\\
&= C\cdot\frac{\Theta(2r)}{\Theta(r)}\cdot\frac{1}{\Theta(2r)}
\int_{B(x_0,2r)}\Phi\left(\frac{|f(x)|}{\sigma}\right)dx\\
&\leq C\cdot\sup_{r>0;B(x_0,r)}\left\{\frac{1}{\Theta(r)}
\int_{B(x_0,r)}\Phi\left(\frac{|f(x)|}{\sigma}\right)dx\right\}.
\end{split}
\end{equation*}
We now turn our attention to the estimate of $K_2$. Recalling that the following estimate holds for any $x\in B$,
\begin{equation*}
\big|\big[b,\mathcal S_\alpha\big](f_2)(x)\big|\leq \big|b(x)-b_{B}\big|\cdot\mathcal S_\alpha(f_2)(x)+\mathcal S_\alpha\Big([b-b_{B}]f_2\Big)(x).
\end{equation*}
Thus, we have
\begin{equation*}
\begin{split}
K_2\leq&\frac{1}{\Theta(r)}\cdot \big|\big\{x\in B:\big|b(x)-b_{B}\big|\cdot\mathcal S_\alpha(f_2)(x)>\sigma/4\big\}\big|\\
&+\frac{1}{\Theta(r)}\cdot \big|\big\{x\in B:\big|\mathcal S_\alpha\Big([b-b_{B}]f_2\Big)(x)\big|>\sigma/4\big\}\big|\\
:=&K_3+K_4.
\end{split}
\end{equation*}
By using the previous pointwise estimate (\ref{S(f2)}), Chebyshev's inequality and the definition of BMO, we conclude that
\begin{equation*}
\begin{split}
K_3&\leq\frac{1}{\Theta(r)}\cdot\frac{\,4\,}{\sigma}\int_B\big|b(x)-b_{B}\big|\cdot\mathcal S_\alpha(f_2)(x)\,dx\\
&\leq C\sum_{j=1}^\infty\frac{1}{|2^{j+1}B|}\int_{2^{j+1}B}\frac{|f(z)|}{\sigma}\,dz
\times\left\{\frac{|B|}{\Theta(r)}\cdot\frac{1}{|B|}\int_B\big|b(x)-b_{B}\big|dx\right\}\\
\end{split}
\end{equation*}
\begin{equation*}
\begin{split}
&\leq C\|b\|_*\sum_{j=1}^\infty\frac{|B|}{|2^{j+1}B|}\cdot\frac{\Theta(2^{j+1}r)}{\Theta(r)}\cdot\frac{1}{\Theta(2^{j+1}r)}\int_{B(x_0,2^{j+1}r)}\frac{|f(z)|}{\sigma}\,dz\\
&\leq C\cdot\sup_{r>0;B(x_0,r)}\left\{\frac{1}{\Theta(r)}
\int_{B(x_0,r)}\Phi\left(\frac{|f(z)|}{\sigma}\right)dz\right\}\times\sum_{j=1}^\infty\frac{|B|}{|2^{j+1}B|}\cdot\frac{\Theta(2^{j+1}r)}{\Theta(r)}.
\end{split}
\end{equation*}
Note that $1\leq D(\Theta)<2^n$, then by using the doubling condition (\ref{doubling}) of $\Theta$, we know that
\begin{equation}\label{theta<C}
\sum_{j=1}^\infty\frac{|B|}{|2^{j+1}B|}\cdot\frac{\Theta(2^{j+1}r)}{\Theta(r)}\leq\sum_{j=1}^\infty\left(\frac{D(\Theta)}{2^n}\right)^{j+1}\leq C,
\end{equation}
which in turn gives that
\begin{equation*}
K_3\leq C\cdot\sup_{r>0;B(x_0,r)}\left\{\frac{1}{\Theta(r)}
\int_{B(x_0,r)}\Phi\left(\frac{|f(z)|}{\sigma}\right)dz\right\}.
\end{equation*}
Applying the previous pointwise estimate (\ref{[b,S](f2)}) and Chebyshev's inequality, we have
\begin{equation*}
\begin{split}
K_4&\leq\frac{1}{\Theta(r)}\cdot\frac{\,4\,}{\sigma}\int_B\Big|\mathcal S_\alpha\Big([b-b_{B}]f_2\Big)(x)\Big|\,dx\\
&\leq\frac{|B|}{\Theta(r)}\cdot\frac{C}{\sigma}\sum_{j=1}^\infty\frac{1}{|2^{j+1}B|}\int_{2^{j+1}B}
\big|b(z)-b_{B}\big|\cdot\big|f(z)\big|\,dz\\
&\leq\frac{|B|}{\Theta(r)}\cdot\frac{C}{\sigma}\sum_{j=1}^\infty\frac{1}{|2^{j+1}B|}\int_{2^{j+1}B}
\big|b(z)-b_{2^{j+1}B}\big|\cdot\big|f(z)\big|\,dz\\
&+\frac{|B|}{\Theta(r)}\cdot\frac{C}{\sigma}\sum_{j=1}^\infty\frac{1}{|2^{j+1}B|}\int_{2^{j+1}B}
\big|b_{2^{j+1}B}-b_B\big|\cdot\big|f(z)\big|\,dz\\
&:=K_5+K_6.
\end{split}
\end{equation*}
For the term $K_5$, notice that the inequality $\Phi(a\cdot b)\leq\Phi(a)\cdot\Phi(b)$ holds for any $a,b>0$, when $\Phi(t)=t(1+\log^+t)$. We then use the generalized H\"older's inequality (\ref{holder}), (\ref{exp}) and (\ref{equiv norm}) together with (\ref{theta<C}) to obtain
\begin{equation*}
\begin{split}
K_5&\leq \frac{|B|}{\Theta(r)}\cdot\frac{C}{\sigma}\sum_{j=1}^\infty\big\|b-b_{2^{j+1}B}\big\|_{\exp L,2^{j+1}B}\big\|f\big\|_{L\log L,2^{j+1}B}\\
&\leq \frac{C\|b\|_*}{\sigma}\cdot\frac{|B|}{\Theta(r)}\sum_{j=1}^\infty
\inf_{\eta>0}\left\{\eta+\frac{\eta}{|2^{j+1}B|}\int_{2^{j+1}B}\Phi\left(\frac{|f(z)|}{\eta}\right)dz\right\}\\
&\leq \frac{C\|b\|_*}{\sigma}\cdot \frac{|B|}{\Theta(r)}\sum_{j=1}^\infty
\left\{\frac{\sigma\cdot\Theta(2^{j+1}r)}{|2^{j+1}B|}+\frac{\sigma}{|2^{j+1}B|}\int_{2^{j+1}B}\Phi\left(\frac{|f(z)|}{\sigma}\right)dz\right\}\\
&\leq C\|b\|_*\cdot\left[1+\sup_{r>0;B(x_0,r)}\left\{\frac{1}{\Theta(r)}
\int_{B(x_0,r)}\Phi\left(\frac{|f(z)|}{\sigma}\right)dz\right\}\right]\\
\end{split}
\end{equation*}
\begin{equation*}
\begin{split}
&\times\sum_{j=1}^\infty\frac{|B|}{|2^{j+1}B|}\cdot\frac{\Theta(2^{j+1}r)}{\Theta(r)}\\
&\leq C\cdot\sup_{r>0;B(x_0,r)}\left\{\frac{1}{\Theta(r)}
\int_{B(x_0,r)}\Phi\left(\frac{|f(z)|}{\sigma}\right)dz\right\}.
\end{split}
\end{equation*}
For the last term $K_6$, an application of the inequality (\ref{BMO}) leads to that
\begin{equation*}
\begin{split}
K_6&\leq C\cdot\frac{|B|}{\Theta(r)}\sum_{j=1}^\infty(j+1)\|b\|_*\cdot\frac{1}{|2^{j+1}B|}\int_{2^{j+1}B}\frac{|f(z)|}{\sigma}\,dz\\
&\leq C\cdot\frac{|B|}{\Theta(r)}\sum_{j=1}^\infty(j+1)\|b\|_*\cdot\frac{\Theta(2^{j+1}r)}{|2^{j+1}B|}\cdot\frac{1}{\Theta(2^{j+1}r)}\int_{B(x_0,2^{j+1}r)}\frac{|f(z)|}{\sigma}\,dz\\
&\leq C\cdot\sup_{r>0;B(x_0,r)}\left\{\frac{1}{\Theta(r)}
\int_{B(x_0,r)}\Phi\left(\frac{|f(z)|}{\sigma}\right)dz\right\}\\
&\times\sum_{j=1}^\infty(j+1)\cdot\frac{|B|}{|2^{j+1}B|}\cdot\frac{\Theta(2^{j+1}r)}{\Theta(r)}.
\end{split}
\end{equation*}
Moreover, by using the doubling condition (\ref{doubling}) of $\Theta$ again and the fact that $1\leq D(\Theta)<2^n$, we find that
\begin{equation}\label{final}
\sum_{j=1}^\infty(j+1)\cdot\frac{|B|}{|2^{j+1}B|}\cdot\frac{\Theta(2^{j+1}r)}{\Theta(r)}\leq C\sum_{j=1}^\infty(j+1)\cdot\left(\frac{D(\Theta)}{2^n}\right)^{j+1}\leq C.
\end{equation}
Substituting the above inequality (\ref{final}) into the term $K_6$, we thus obtain
\begin{equation*}
K_6\leq C\cdot\sup_{r>0;B(x_0,r)}\left\{\frac{1}{\Theta(r)}
\int_{B(x_0,r)}\Phi\left(\frac{|f(z)|}{\sigma}\right)dz\right\}.
\end{equation*}
Summing up all the above estimates, we therefore conclude the proof of the main theorem.
\end{proof}

Finally, we remark that by using the same arguments as in the proof of Theorems \ref{mainthm:6} and \ref{mainthm:7}, we can also show the conclusion of Theorem \ref{mainthm:9}. The details are omitted here.

\section{Acknowledgment}

This paper is supported by the Fundamental Research Funds for the Central Universities of Hunan University (Grant no. 531107040013).

\end{document}